%% file: main_arxiv.tex
\documentclass[11pt]{article}
\usepackage{setspace}
\usepackage{multicol}
\usepackage{adjustbox}
\usepackage{pdflscape}
\usepackage{rotfloat}
\usepackage{makecell}
\usepackage{diagbox}
\usepackage{graphicx}
\usepackage[a4paper, left=1in, right=1in, top=1in, bottom=0.8in]{geometry}
\usepackage{longtable}

\usepackage{amssymb}
\usepackage{amsmath}
\usepackage{amsthm}

\usepackage{natbib}
 \bibpunct[, ]{(}{)}{,}{a}{}{,}%

\usepackage{algorithm}
\usepackage{calc}
\usepackage{todonotes}
\usepackage[noend]{algpseudocode}
\allowdisplaybreaks
\usepackage[colorlinks=false,linktoc=section]{hyperref}
\usepackage{array}
\usepackage{xspace}
\usepackage{mathtools}
\usepackage{subfig}
\usepackage{booktabs}
\usepackage{cleveref}

\newcommand{\LineComment}[1]{\hfill $\triangleright$ #1}
\newcommand{\sep}{,\ }

\begin{document}

\title{An Adaptive Variable Neighborhood Search for a Family of Set Covering Routing Problems with an Application in Disaster Relief Operations}

\author{
  Andreas Hagn\thanks{Corresponding author. Technical University of Munich, Chair of Operations Management, Arcisstr.~21, 80333 München, Germany. \texttt{andreas.hagn@tum.de}} \and
  Jan Krause\thanks{University of Technology Nuremberg, Analytics and Optimization, Dr.-Luise-Herzberg-Straße~4, 90461 Nuremberg. \texttt{jan.krause@utn.de}} \and
  Moritz Stargalla\thanks{University of Technology Nuremberg, Applied Discrete Mathematics, Dr.-Luise-Herzberg-Straße~4, 90461 Nuremberg. \texttt{moritz.stargalla@utn.de}} \and
  Lorenza Moreno\thanks{Federal University of Juiz de Fora, Computer Science Department, Rua José Lourenço Kelmer, s/n São Pedro, Juiz de Fora - MG, 36036-900 Brazil. \texttt{lorenza.moreno@ufjf.br}}
}

\date{}

\maketitle

\begin{abstract}
This paper studies a variant of the Set Covering Routing Problem (SCRP) motivated by post-disaster humanitarian logistics. We consider a hybrid distribution concept in which the majority of transportation is performed by helicopters, while ground transport is limited to the last mile, addressing severe accessibility constraints in disaster-affected regions. The resulting problem integrates landing site location, routing, and covering decisions, incorporating features of the Multi-Vehicle Covering Tour Problem (m-CTP) and the Vehicle Routing with Demand Allocation Problem (VRDAP) in a facility-capacitated, multi-depot setting.
Due to the computational complexity of the problem, we develop an Adaptive Variable Neighborhood Search (AVNS) that combines established routing operators with novel mechanisms for covering decisions.
The performance of the proposed approach is evaluated on benchmark instances for the related m-CTP and VRDAP problems, demonstrating competitive solution quality compared to problem-specific state-of-the-art approaches. Furthermore, we apply our AVNS to a real-world case study based on the 2024 flash floods in Afghanistan. The results highlight the practical relevance of the proposed framework and provide managerial insights into effective distribution strategies for disaster response operations.
\end{abstract}


\bigskip

\input{sections/introduction}

\input{sections/literature}

\input{sections/problem_description}

\input{sections/algorithm}

\input{sections/computational_experiments}
\input{sections/case_study}

\input{sections/conclusion}

\input{sections/acknowledgement}

\noindent \textbf{Declaration of generative AI and AI-assisted technologies in the manuscript preparation process} \\
During the preparation of this work the authors used ChatGPT and DeepL for conceptualization, literature study and language refinement. After using these services, the authors reviewed and edited the content as needed and take full responsibility for the content of the publication.

\bibliographystyle{apalike}
\bibliography{bibliography}

\clearpage
\appendix
\input{sections/supplementary}

\end{document}

%% file: sections/introduction.tex
\section{Introduction}\label{sec:introduction}

According to United Nations Office for the Coordination of Humanitarian Affairs \citep{UNOCHA2024GHO2025-july}, around 300 million people were predicted to be in need of humanitarian aid in 2025, whereas only 181.2 million people can be targeted by humanitarian aid operations. At the same time, the World Food Programme (WFP) has predicted an alarming funding shortfall of 34 percent for 2025 compared to 2024 \citep{wfp_funding_shortfall}. While the majority of global humanitarian needs stem from protracted crises driven by political instability, armed conflict, and economic collapse (e.g., \cite{UNOCHA2024GHO2025} for an overview of conflicts), sudden-onset natural disasters trigger acute challenges for humanitarian logistics. They are increasingly linked to climate change and continue to create urgent, localized demands for rapid response and logistics coordination. As with all natural disasters, floods pose the particular challenge of enabling rapid response operations to deliver life-saving relief items to affected populations. A major difficulty in flood emergencies arises from the often widespread inaccessibility of the transport infrastructure, as large parts of the road network can become impassable or even completely destroyed, severely constraining logistical accessibility and the timely distribution of aid. Developing more reliable distribution networks is therefore essential. One promising approach is to cover a substantial part of the transport distance by air, while limiting ground transport to the last mile. More precisely, we consider the following distribution concept for humanitarian relief goods: airports serve as central depots from which helicopters conduct circular flights, delivering bulk consignments of relief items to local landing sites, i.e., distribution points. Beneficiaries then travel from their location to their assigned distribution point and collect their aid, supervised by local volunteers. Our focus lies on the post-flood phase, when the immediate emergency response has passed, and the road network's condition has improved sufficiently to enable light transport by foot. To obtain an efficient solution to this distribution problem, we propose an extended vehicle routing problem that incorporates covering aspects.
\\
For this purpose, we introduce a new model formulation addressing the following problem: 
Given a set of customers (affected individuals), a set of potential depots (airports), a set of potential facilities (landing sites), and a limited fleet of capacitated vehicles (helicopters), the goal is to 
\begin{itemize}
\setlength{\itemsep}{0pt}
\setlength{\parskip}{0pt}
\setlength{\parsep}{0pt}
    \item select a subset of facilities,
    \item assign each customer to exactly one selected facility that can cover them, and
    \item construct vehicle routes that start and end at depots while visiting the selected facilities along the routes,
\end{itemize}
so as to satisfy all customer demands at minimum routing and assignment costs, subject to vehicle and facility capacity constraints as well as route-length limitations.
On the one hand, the considered problem can be classified as a \textit{Set Covering Routing Problem} (SCRP) as defined in \cite{MORADI2024110730}, which can be considered a generalization of the \textit{Multi-Vehicle Covering Tour Problem} (m-CTP) and the \textit{Vehicle Routing with Demand Allocation Problem} (VRDAP). To the best of our knowledge, no existing study in the literature addresses the proposed problem setting.
Furthermore, \cite{MORADI2024110730} emphasizes that most studies on the SCRP focus on single-depot scenarios, likely due to the already high inherent complexity introduced by the covering aspects.
Nevertheless, extending the research to multi-depot settings is a crucial step toward enhancing the practical applicability of routing algorithms in this field. 
To efficiently solve this problem for realistically sized instances, we develop an \textit{Adaptive Variable Neighborhood Search} (AVNS) integrating several well-known mechanisms, such as tabu simulated annealing and adaptive operator selection.
As local search operators, we employ well-established operators from the vehicle routing literature and combine them with newly developed ones designed to improve depot selection and customer-to-facility assignment. 
The adaptive operator selection mechanism, inspired by the remarkable success of \textit{Adaptive Large Neighborhood Search} (ALNS), enhances the robustness of our heuristic across different instance types. 
To assess the potential of our developed AVNS for the considered problem in the context of the initially motivated flood response application, we perform a realistic case study on the 2024 flash floods in Afghanistan. According to \cite{unicef_extreme_weather}, Afghanistan ranks among the ten countries most severely affected by extreme weather conditions and natural disasters, including droughts, storms, avalanches and earthquakes. The flash appeal released by the WFP in \cite{wfp_heb_response} outlines the severe impacts of the heavy rainfall between March and May 2024, which caused extensive flooding across large parts of Afghanistan. Thousands of households and tens of thousands of people were affected, with hundreds reported dead or injured. The magnitude of the situation necessitated an urgent funding appeal of USD 14.5 million to assist 80,000 people. Declining funding levels, however, posed significant threats to the WFP, underscoring the need for efficient management of emergency responses. Multiple reports from several humanitarian organizations document the inaccessibility of the transport network, including damaged roads and bridges, revealing the complexity of ground-level aid distribution (e.g., see \cite{inaccessible_source_who}, \cite{inaccessible_source_unicef}, \cite{awd_cholera}). In addition to the flood-induced deterioration of ground transport conditions, the existence of numerous hard-to-reach communities in mountainous regions and persistent security concerns due to political instability further highlight the relevance of supporting humanitarian logistics through air transport in Afghanistan. Thus, we apply the proposed algorithm to the aforementioned application. In this context, we derive managerial insights into optimal fleet sizes, trade-offs between operational costs, last-mile burdens and fairness, and investigate the potential of concepts commonly used in humanitarian aid to reduce costs.\\
In summary, our methodological and application-oriented contributions to humanitarian aid distribution in the context of flood response can be outlined as follows: 
\begin{itemize}
\setlength{\itemsep}{0pt}
\setlength{\parskip}{0pt}
\setlength{\parsep}{0pt}
    \item[i)] Inspired by a flood response application, we introduce a novel multi-depot generalization of two known variants of SCRPs, which integrates facility selection, customer covering, and capacitated vehicle routing decisions in a unified approach.
    \item[ii)] To efficiently solve this complex problem, we design an AVNS that integrates Tabu Search and Simulated Annealing mechanisms, combining well-established local search operators with newly developed ones.
    \item[iii)] Benchmarks of the proposed heuristic demonstrate its competitiveness relative to problem-specific state-of-the-art approaches of both the m-CTP and the VRDAP.
    \item[iv)] We demonstrate the potential of our approach through an extensive case study on the 2024 flash floods in Afghanistan, providing substantial insights for practitioners.
\end{itemize}

The remainder of the paper is organized as follows. \Cref{sec:literature_new} positions our work within the literature from modeling, methodological, and application perspectives.
\Cref{sec:problem_description} provides a formal problem description and its formulation as a mixed-integer linear program (MILP). \Cref{sec:algorithm} presents the proposed AVNS in detail, while \Cref{sec:comp_experiments} reports results obtained from applying the proposed approach to several benchmark instance sets. In \Cref{sec:case_study}, we present the aforementioned case study and derive several insights into the potential and limitations of the proposed distributional structure. Finally, in \Cref{sec:conclusion}, we summarize our findings and provide potential future research directions.

%% file: sections/literature.tex
\section{Literature Review}
\label{sec:literature_new}

In the following, we provide a thorough overview of relevant literature. In \Cref{sec:scrp_problems}, we position the problem at hand in the class of SCRPs and briefly address their core features. Because the m-CTP and the VRDAP can be considered special cases of the problem addressed in this publication, \Cref{sec:literature_mctp} and \Cref{sec:literature_vrdap} summarize key publications for said problems. In \Cref{sec:avns_survey}, we showcase that VNS-based heuristics provide competitive performance for vehicle routing problems and its variants. Finally, \Cref{sec:hum_aid_routing} provides an overview of core contributions applying vehicle routing methods to problem settings arising in humanitarian aid.

\subsection{Set Covering Routing Problems}\label{sec:scrp_problems}
Problems combining vehicle routing with covering aspects have seen an increase in interest in the past years. Due to the sheer volume of existing literature we follow the classification scheme proposed in \cite{MORADI2024110730}. Their studied problem class comprises a main subclass of so-called non-hamiltonian cyclic routing problems, i.e., not every customer has to be visited directly, and routes start and end at the same depot. There are different handling strategies with respect to the non-visited customers. First, they can be left isolated or assigned a profit/prize to be on the routes which leads to the subclass of routing problems with profits. Second, they can either be visited, covered by intermediary nodes (called \textit{facilities}) or by other customer nodes, or left isolated (Maximal Covering Routing Problems). Last, SCRPs correspond to the class of problems where the non-visited customers may be directly visited or covered and no isolation is feasible. The authors present a unified notation and general mathematical model with the aim of covering all SCRP variants. Based on the defined classification scheme their study includes reviews of ten SCRP variants that the authors consider to be most-studied and latest in prevalent literature.
Notably, their classification highlights a limitation of the existing literature: multi-depot problems remain largely underexplored, despite their practical relevance in applications such as humanitarian logistics. To be precise, only  \cite{allahyari2015hybrid} and \cite{nedjati2017bi} deal with multi-depot SCRPs.\\
\begin{table}[h!]
\centering
\begin{tabular}{lccc}
\hline
Feature & m-CTP & VRDAP & our problem \\
\hline
multi-depot           &   &   & X \\
vehicle capacity          &   & X & X \\
facility capacity           &   &  & X \\
coverage restrictions      & X  &   & X \\
facility limit per route   & X  &  & X \\
distance limit per route        & X &  & X \\
route cost minimization        & X & X & X \\
assignment cost minimization  &   & X & X \\
\hline
\end{tabular}
\caption{Comparison of problem features between m-CTP, VRDAP and our problem}
\label{tab:problem_comparison}
\end{table}
While many SCRPs share several similarities, two of its variants are special cases of the problem considered in this publication: the m-CTP and the VRDAP. A comparison of problem features between our considered problem and the aforementioned problems can be found in \Cref{tab:problem_comparison}. For a more comprehensive comparison of SCRPs and their relationships, the interested reader is referred to Section~3.13 in \cite{MORADI2024110730}. In the following two sections, we provide a literature review of these two problem types. For an overview of practical applications of SCRPs, we refer to the publications referenced in \cite{MORADI2024110730}. Furthermore, SCRPs have gained significant traction in the context of out-of-home delivery, see, e.g., \cite{janinhoff2024out}.

\subsection{Multi-Vehicle Covering Tour Problem}\label{sec:literature_mctp}
The \textit{Covering Tour Problem} (CTP) was originally introduced by \cite{Gendreau1997} and aims to find an assignment of customers to facilities, as well as a route that visits all selected facilities, such that the total routing cost is minimized. In this context, for each customer, there exists a subset of facilities to which they can be assigned. In contrast to the problem considered in this publication, facilities, as well as vehicles are assumed to be uncapacitated and covering does not incur any cost. A subset of facilities can be mandatory, i.e., they have to be visited even if no customer is assigned to them. \cite{Hachicha2000} are the first to consider a multi-vehicle variant of this problem, further enhancing the practical applicability of the problem by also introducing upper bounds on the maximum number of stops and the maximum distance per route. 
Note that the former condition can be interpreted as an implicit vehicle capacity, assuming unit demands for customers.
For a comprehensive literature review, the interested reader is referred to \cite{Oliveira2025_mCTP_exact}. Note that several papers address a variant of this problem without restrictions on maximum route distances, which is typically abbreviated as m-CTP-p. While  \cite{Glize2020}, \cite{Jozefowiez2014}, and \cite{Oliveira2025_mCTP_exact} solve the m-CTP and the m-CTP-p using exact solution methods, \cite{Kammoun2016} and \cite{Oliveira2015} propose heuristic solution. 

\subsection{Vehicle Routing with Demand Allocation Problem}\label{sec:literature_vrdap}
The \textit{Vehicle Routing with Demand Allocation Problem} (VRDAP) was initially introduced in \cite{Ghoniem2013} motivated by logistical operations for a food bank. In contrast to the m-CTP, the VRDAP does not impose any restriction on the maximum number of nodes or the maximum distance per route. Furthermore, vehicles are assumed to be capacitated and covering decisions incur costs. Moreover, customer-to-facility assignments are unrestricted, i.e., each customer can be assigned to each facility. While there exist several related problems in the literature (cf. \cite{MORADI2024110730}), the number of publications directly adressing VRDAP is rather limited. \cite{Oliveira2025_VRDAP_exact} and \cite{Reihaneh2019} solve the problem exactly, whereas heuristic approaches are described in \cite{Ghoniem2013}, \cite{Oliveira2025_VRDAP_ILS}, and \cite{Reihaneh2018}.\\
Overall, neither the m-CTP nor the VRDAP captures the full set of decision dimensions required in our application. This gap motivates the unified problem setting considered in this paper.

\subsection{Heuristics in Related Problems}\label{sec:avns_survey}
Over the past decades, a variety of metaheuristic frameworks have been established for solving numerous variants of vehicle routing problems, including single-solution–based methods such as \textit{Iterated Local Search} and \textit{(Adaptive) Large Neighborhood Search} ((A)LNS), as well as population-based approaches such as \textit{Genetic Algorithms} and \textit{Memetic Algorithms}.
While it is often difficult to identify the best performing heuristic framework for a variant of a vehicle routing problem without conceptualizing and implementing all of them, \cite{heuristicssurvey} surveys more than 300 publications on heuristics in vehicle routing problems. They highlight the trend of combining various heuristic concepts rather than strictly following a single one. \cite{Vidal2013} come to a similar conclusion. Moreover, \cite{MORADI2024110730} showcase in Section~3.14 that methods based on \textit{Variable Neighborhood Search} (VNS), which was introduced in \cite{Mladenovi1997}, frequently rank among the most competitive frameworks for SCRPs. For these reasons, we utilize a general VNS framework and enhance it by several concepts that are usually applied to ALNS methods. For a comprehensive introduction as well as recent trends in the area of VNS, the interested reader is referred to \cite{bookmetaheuristics2019} and \cite{Brimberg2023}.

\subsection{Routing in Humanitarian Aid}\label{sec:hum_aid_routing}
 For handling the discrepancy of high humanitarian needs on the one hand with comparably few resources on the other hand, increasing efficiency in relief operations is of pivotal importance. The latter can be well adressed by OR methods, making humanitarian aid an acitve branch in OR literature. A multitude of literature surveys highlighting different aspects of modeling or solution techniques exists, e.g., \cite{AnayaArenas2014}, \cite{Anuar2021}, \cite{Hoyos2015}, \cite{Luo2023}, \cite{delaTorre2012}, \cite{Xu2018}, \cite{Zheng2015}. 
It can be observed that the majority of OR literature in humanitarian aid focuses on acute interventions due to emergency events rather than on post-acute relief operations. 
 In particular, vehicle routing problems play a central role in the humanitarian supply chain. 
More precisely, the covering tour concept is inevitable for routing in disaster areas due to damaged road infrastructure and the need for the timely provision of relief goods to affected people. Some of the publications including covering aspects are \cite{Alinaghian2019}, \cite{Davoodi2019}, \cite{Goli2019}, \cite{Kl2025}, \cite{Nolz2010}, and \cite{Tricoire2012}. In line with the latter references, structurally different objective functions are used in routing models for humanitarian aid. In highly time-critical settings, the focus is typically on minimizing the arrival time of the last beneficiary or the total latency (see also \cite{MoshrefJavadi2016}). In contrast, in post-acute scenarios such as the case study considered in this paper, it is often important to reduce distribution costs and improve aid accessibility, for instance by minimizing the total transportation distance, including last-mile distances.

%% file: sections/problem_description.tex
\section{Problem Description}
\label{sec:problem_description}

\noindent 
In this section, we first provide a formal introduction of the considered problem and subsequently present a mixed-integer linear programming formulation. Using established terminology from the literature, our problem can be described as a multi-depot, facility-capacitated version of the VRDAP incorporating route-length limitations and coverage restrictions as in the \mbox{m-CTP}. Formally, let $D$ denote the set of available depots and $F$ the set of potential facilities (delivery sites), from which customers can pick up their goods. Based on these sets, we define a set of arcs as
\[
A \coloneqq \{(i, j) : i \in D \cup F,\, j \in F,\ i\neq j\} \cup \{(i, j) : i \in F,\, j \in D \cup F, i \neq j\}.
\]
Furthermore, let $W$ be the set of customers, each with a non-negative demand $r_l \geq 0$ for all $l\in W$ of a homogeneous good. We assume that each customer $l \in W$ can be assigned only to a subset of facilities $\emptyset \neq F_l \subseteq F$. Each facility $i \in F$ can serve at most a total demand of $Q^F_i \in (0,\infty]$. Moreover, a fleet of up to $M$ vehicles with a homogeneous maximum capacity $Q \in (0,\infty]$ is available. Beyond these capacity constraints, each route is subject to an upper bound on the number of facilities $p\in\mathbb{N}$ and on the total travel distance $q\geq 0$, where the distance of an arc $(i,j)\in A$ is denoted by $\text{dist}_{ij} \geq 0$.
In addition, routing costs $c_{ij} \ge 0$ are incurred for traversing arc $(i, j) \in A$, and assignment costs $a_{li} \ge 0$ arise when customer $l \in W$ is assigned to facility $i \in F_l$.
Then, the problem consists of selecting a subset $\mathcal{F}\subseteq F$ of facilities, a surjective assignment function $\phi: W\rightarrow \mathcal{F}$ of customers to selected facilities and constructing up to $M$ vehicle routes $\mathcal{R}$, one for each vehicle, where each route in $\mathcal{R}$ starts at a depot in $\mathcal{D}$, visits facilities in $\mathcal{F}$ and returns back to its starting depot, such that 
{
\begin{itemize}
\setlength{\itemsep}{0pt}
\setlength{\parskip}{0pt}
\setlength{\parsep}{0pt}
    \item each customer $l\in W$ is assigned to exactly one facility in $F_l$,
    \item each facility in $\mathcal{F}$ is visited by exactly one vehicle,
    \item each route in $\mathcal{R}$ does not exceed vehicle capacity $Q$, maximum facility number $p$ and maximum total distance $q$,
    \item the demand of customers assigned to each facility $i\!\in\!\mathcal{F}$ does not exceed its capacity $Q_i^F$, and
    \item the total costs, i.e., the sum of routing and assignment costs, are minimized.
\end{itemize}
}

\Cref{fig:problem_example} illustrates a feasible solution to the considered problem with two depots, seven facilities visited by three vehicle routes, three additional non-routed facilities, and assignments of the 19 customers to the visited facilities. Note that the sets $F_l$ containing the covering possibilities for customers $l \in W$ are implicitely given by covering radii. 

\begin{figure}[H]
    \centering
    \begin{tikzpicture}[scale=0.7, every node/.style={font=\small}]

\definecolor{depotcol}{RGB}{227,170,180}
\definecolor{facilitycol}{RGB}{84,158,104}
\definecolor{customercol}{RGB}{126,98,168}
\definecolor{covercol}{RGB}{160,170,210}

\tikzset{
    depot/.style={
        circle, draw=black, fill=depotcol,
        minimum size=9mm, inner sep=0pt, font=\bfseries
    },
    facility/.style={
        rectangle, draw=black, fill=facilitycol,
        minimum size=5.5mm, inner sep=0pt, font=\scriptsize
    },
    facilityNR/.style={
        rectangle, draw=black!50, fill=gray!35,
        text=black!60,
        minimum size=5.5mm, inner sep=0pt, font=\scriptsize
    },
    customer/.style={
        circle, draw=black, fill=customercol,
        minimum size=3.7mm, inner sep=0pt, font=\scriptsize, text=white
    },
    routeA/.style={line width=1.0pt},
    routeB/.style={line width=1.0pt},
    assign/.style={dotted, line width=0.8pt},
    cover/.style={
    draw=covercol!50,
    fill=covercol,
    fill opacity=0.03,
    line width=0.3pt
}
}

\node[depot] (D1) at (0.3,4.2) {D};
\node[depot] (D2) at (12.8,4.2) {D};

\node[facility] (F1) at (2.2,6.3) {8};
\node[facility] (F2) at (3.8,4.4) {9};
\node[facility] (F3) at (2.4,2.1) {7};

\node[facility] (F4) at (7.1,6.5) {8};
\node[facility] (F5) at (8.9,4.5) {9};
\node[facility] (F6) at (7.1,2.2) {6};

\node[facility] (F7) at (10.8,2.2) {5};

\node[facilityNR] (F8)  at (5.5,7.8) {6};
\node[facilityNR] (F9)  at (5.6,0.9) {5};
\node[facilityNR] (F10) at (10.6,6.9) {4};

\draw[cover] (F1) circle (1.45);
\draw[cover] (F2) circle (1.45);
\draw[cover] (F3) circle (1.45);

\draw[cover] (F4) circle (1.45);
\draw[cover] (F5) circle (1.45);
\draw[cover] (F6) circle (1.45);

\draw[cover] (F7) circle (1.35);

\draw[cover] (F8) circle (1.25);
\draw[cover] (F9) circle (1.25);
\draw[cover] (F10) circle (1.25);

\node[customer] (c1)  at (1.3,7) {2};
\node[customer] (c2)  at (2.2,7.5) {3};
\node[customer] (c3)  at (3.1,7.0) {2};

\node[customer] (c4)  at (3.0,5) {4};
\node[customer] (c5)  at (4.6,4.8) {2};
\node[customer] (c6)  at (3.9,3.4) {2};

\node[customer] (c7)  at (1.5,1.6) {2};
\node[customer] (c8)  at (2.4,0.9) {3};
\node[customer] (c9)  at (3.2,1.7) {1};

\node[customer] (c10) at (6.1,7) {2};
\node[customer] (c11) at (7.2,7.7) {3};
\node[customer] (c12) at (8,7.1) {2};

\node[customer] (c13) at (8.0,5) {3};
\node[customer] (c14) at (9.8,5.2) {2};
\node[customer] (c15) at (9.2,3.5) {3};

\node[customer] (c16) at (6.2,1.7) {2};
\node[customer] (c17) at (7.5,1.2) {3};

\node[customer] (c18) at (10.0, 2.0) {2};
\node[customer] (c19) at (11.5, 2.0) {2};

\draw[assign] (c1) -- (F1);
\draw[assign] (c2) -- (F1);
\draw[assign] (c3) -- (F1);

\draw[assign] (c4) -- (F2);
\draw[assign] (c5) -- (F2);
\draw[assign] (c6) -- (F2);

\draw[assign] (c7) -- (F3);
\draw[assign] (c8) -- (F3);
\draw[assign] (c9) -- (F3);

\draw[assign] (c10) -- (F4);
\draw[assign] (c11) -- (F4);
\draw[assign] (c12) -- (F4);

\draw[assign] (c13) -- (F5);
\draw[assign] (c14) -- (F5);
\draw[assign] (c15) -- (F5);

\draw[assign] (c16) -- (F6);
\draw[assign] (c17) -- (F6);

\draw[assign] (c18) -- (F7);
\draw[assign] (c19) -- (F7);

\draw[routeA] (D1) -- (F1) -- (F2) -- (F3) -- (D1);
\draw[routeB] (D2) -- (F4) -- (F5) -- (F6) -- (D2);
\draw[routeA] (D2) -- (F7) -- (D2);

\end{tikzpicture}
    \caption{Illustrative example of the considered problem with depots (red), capacitated routed and non-routed facilities (green and grey), customers with demands (purple) and customer covering radii (light blue shaded circles)}
    \label{fig:problem_example}
\end{figure}

For the remainder of the paper, a feasible solution to the considered problem is denoted as $(\mathcal{R}, \phi)$. We say that a customer $l\in W$ is \textit{covered} by a facility $i\in F_l$ if $\phi(l) = i$ holds. Furthermore, the sets $\mathcal{D}$ and $\mathcal{F}$ are called the sets of \textit{selected depots} and \textit{selected facilities}, respectively. Note that we do not explicitly model facility opening costs for facilities $i \in F$, as these can be modeled as part of the arc traversal costs $c_{ij}$ for traversing to the subsequent depot or facility $j \in D \cup F$. 
Note that mandatory facilities and direct customer visits can be modeled indirectly using dummy customers or facilities, respectively.

\subsection{A Mixed-Integer Linear Programming Formulation}
\label{subsec:MIP}
In the following, we formulate the problem as described in \Cref{sec:problem_description} as a mixed-integer linear program. The proposed formulation is essentially derived from the unified mathematical model for SCRPs presented in \cite{MORADI2024110730}. For all arcs $(i, j) \in A$, variables $x_{ij}\in\{0,1\}$ indicate whether they are traversed by a vehicle or not. Furthermore, variables $z_{li}\in\{0,1\}$ encompass the allocation of customers $l \in W$ to facilities $i \in F_l$. In addition, $w_i\in\{0,1\}$ denotes whether facilities $i \in F$ are open. Finally, continuous variables $u_i\geq0$, $\pi_i\geq0$, and $\delta_i\geq0$ are introduced to track vehicle loads, and to model the $p$- and $q$-constraints. The model can now be formulated as follows: 

\begin{align}
\min &
 \sum_{(i, j) \in A} c_{ij}\,x_{ij}
+ \sum_{\substack{l \in W \\ i \in F_l}} a_{li}\,z_{li} & \label{MIP:obj}\\
\text{s.t. } & \sum_{d\in D}\sum_{i\in F} x_{di} \le M, \label{MIP:max-vehicle1}\\
& \sum_{d\in D}\sum_{i\in F} x_{id} \le M, \label{MIP:max-vehicle2}\\
& \sum_{i\in F} x_{id}
  = \sum_{i\in F} x_{di},
  & \forall d\in D, \label{MIP:routing1}\\
& \sum_{ (i, j) \in A } x_{ij}
  = \sum_{ (j, i) \in A } x_{ji}
  = w_i, & \forall i\in F, \label{MIP:routing2}\\
  & \sum_{i \in F_l}
 z_{li} = 1,
  & \forall l \in W, \label{MIP:covering}\\
& z_{li} \le w_i,& \forall l \in W,\;
  \forall i \in F_l, \label{MIP:zw-link}\\
& \sum_{l\in W: i \in F_l} r_l\, z_{li} \le Q^F_i,
  & \forall i\in F, \label{MIP:facil-cap}\\
  & u_{j}  \le u_{i} - 
\sum_{l\in W: i \in F_l} r_l\, z_{li} + Q(1 - x_{ij}), &
  \forall i,j\in F,\; i\neq j, \label{MIP:vehicle-cap1}\\
  & \sum_{l \in W: i \in F_l}r_l z_{li} \le u_i \le Q, & \forall i \in F, \label{MIP:vehicle-cap2}\\
& \pi_{j} \ge \pi_{i} + x_{ij} - (p-1)(1 - x_{ij}),&
  \forall i\in F,\;
  j\in F,\; i\neq j, \label{MIP:p-constr1}\\
& \pi_{i} \le p - 1, &
  \forall i\in F, \label{MIP:p-constr2}\\
  & \delta_{j} \ge \delta_{i} + \text{dist}_{ij}\, x_{ij} - q(1 - x_{ij}), &
  \forall i\in F,\;
  j\in F,\; i\neq j, \label{MIP:q-constr1}\\
  & \delta_{i} \ge \text{dist}_{di}\, x_{di}, & \forall i\in F,\; \forall d\in D, \label{MIP:q-constr2}\\
& \delta_{i} + \text{dist}_{id}\, x_{id} \le q, &
  \forall i\in F,\; \forall d\in D, \label{MIP:q-constr3}\\
& x_{ij} \in \{0,1\},
& \forall (i,j)\in A, \label{MIP:x-var}\\
& z_{li} \in \{0,1\}, & \forall l \in W, \; \forall i \in F_l, \label{MIP:z-var}\\
& w_i \in \{0,1\}, & \forall i \in F, \label{MIP:w-var}\\
& u_i \ge 0, & \forall i \in F, \label{MIP:u-var}\\
& \pi_i \ge 0, & \forall i \in F,\label{MIP:pi-var} \\
& \delta_i \ge 0, & \forall i \in F. \label{MIP:delta-var}
\end{align}

Line \eqref{MIP:obj} contains the objective function, which minimizes the sum of routing costs and the costs for assigning customers to facilities. Constraints \eqref{MIP:max-vehicle1} and \eqref{MIP:max-vehicle2} limit the number of vehicles used to at most $M$. Constraints \eqref{MIP:routing1} and \eqref{MIP:routing2} are flow conservation constraints. Constraints \eqref{MIP:covering} ensure that each customer is covered by exactly one facility. Constraints \eqref{MIP:zw-link} link the $z$- and $w$-variables, requiring that customers can only be assigned to opened facilities, and facilities must be opened as soon as at least one customer is assigned to them. Constraints \eqref{MIP:facil-cap} describe the facility capacities. Constraints \eqref{MIP:vehicle-cap1} and \eqref{MIP:vehicle-cap2} model vehicle capacities, which implicitly also eliminate subtours. Here, for each $i\in F$, $u_i$ can be interpreted as the (leftover) vehicle load before serving $i$. Note that constraints \eqref{MIP:vehicle-cap1} become redundant if arc $(i, j)$ is not traversed. Similarly, constraints \eqref{MIP:p-constr1} and \eqref{MIP:p-constr2}, as well as constraints \eqref{MIP:q-constr1}-\eqref{MIP:q-constr3}, represent the $p$- and $q$-constraints, respectively. In a feasible solution, $\pi_i$ and $\delta_i$ can be interpreted, for all facilities $i \in F$, as the number of already visited facilities, including $i$, and the already traveled distance, respectively. Again, constraints \eqref{MIP:p-constr1} and \eqref{MIP:q-constr1} become redundant when arc $(i, j)$ is not traversed. Finally, constraints \eqref{MIP:x-var}–\eqref{MIP:delta-var} define the decision variables of the optimization problem.\\
While the proposed model can be solved for small-sized instances in reasonable time, it becomes computationally intractable for instances with several hundred customer and facility nodes, highlighting the necessity for an efficient heuristic solution method.

%% file: sections/algorithm.tex
\section{Algorithm}
\label{sec:algorithm}
In this section, we present an adaptive variable neighborhood search to solve the problem at hand. Since its introduction in \cite{Mladenovi1997}, VNS has become one of the most successful metaheuristics, with numerous applications in combinatorial and integer optimization. In particular, the AVNS proposed in the following can be seen as an enhancement of the basic concept of a VNS. In general, a VNS iteratively improves a given initial solution over time through a systematic alternation of exploration and exploitation. In the exploratory phase, also referred to as perturbation or shaking, the main goal is to move into different regions of the solution space in order to escape local optima. In contrast, the exploitation phase applies local search operators to further enhance the solution. Both mechanisms are implemented within distinct neighborhood structures. For shaking, this implies that different intensities of exploration can be realized, while for the local search, various operators are available that can lead to improved solutions in different ways. 
Our algorithm extends the VNS by incorporating several components that are largely inspired by the class of ALNS, e.g., see \cite{WindrasMara2022}. To be precise, we employ a smoothed scoring mechanism, as well as a simulated annealing (SA) step, similar to \cite{Ropke2006} and \cite{Sacramento2019}, which randomizes operator selection in each iteration and biases it toward selecting operators that have frequently yielded good solutions in past iterations. Furthermore, we utilize a tabu list to detect and break cyclic operator behaviour.\\
In \Cref{sec:overall_framework}, the proposed AVNS framework is presented. In \Cref{sec:constructive_heuristic} throughout \Cref{sec:shaking}, we provide a closer examination of the elements that are central to the proposed algorithm. In \Cref{sec:local_search}, we present the operators used in the proposed algorithm. This includes a novel string replacement operator, which iteratively attempts to remove facility strings and reassign customers to open or closed facilities.

\subsection{AVNS Framework}\label{sec:overall_framework}
Initially, the AVNS constructs a feasible solution derived using the heuristic presented in \Cref{sec:constructive_heuristic}, and then iteratively attempts to improve it. \Cref{alg:avns_framework} summarizes the core mechanism of the proposed AVNS. 

\begin{algorithm}[H]
\caption{Adaptive Variable Neighborhood Search (AVNS)}\label{alg:avns_framework}
\begin{algorithmic}[1]
\Require problem instance $\mathcal{I}$, local search operators $\mathcal{O}$,
         initial score $\lambda^{\text{initial}} > 0$,
         score updates $\sigma \in \mathbb{R}_{\geq 0}^4$,
         initial temperature $T^{\text{initial}} > 0$,
         constructive heuristic calls $N^{\text{start}}$,
         smoothing factor $\rho \in [0,1]$,
         segment length $N^{\text{seg}} \in \mathbb{N}$,
         tabu list size $N^{\text{tabu}} \in \mathbb{N}$,
         stagnation limit $N^{\text{stag}} \in \mathbb{N}$,
         shaking attempts $N^{\text{shaking}}\in\mathbb{N}$,
         shaking neighborhoods $(\Lambda_k)_{k=1,\ldots,K}$,
         maximum number of checks $n_{\text{max}}$ per facility within FacilityStringReplacement,
         maximum runtime $t^{\max} > 0$
\State $(\mathcal{R}, \phi) \gets \text{ConstructiveHeuristic}(\mathcal{I}, N^{\text{start}})$
\If{no initial solution found}
    \State \Return None
\EndIf
\State $(\mathcal{R}^{\text{best}}, \phi^{\text{best}}) \gets (\mathcal{R}, \phi)$
\State $T, T^{\text{start}} \gets T^{\text{initial}} \cdot \text{cost}(\mathcal{R}, \phi)$ \Comment{start temperature for tabu simulated annealing}
\State $\mathcal{L} \gets [(\mathcal{R}, \phi)]$ \Comment{tabu list}
\State $\lambda \gets (\lambda^{\text{initial}})_{o \in \mathcal{O}}, \lambda^{\text{seg}} \gets (0)_{o\in\mathcal{O}}$ \Comment{scores and segment scores}
\State $k \gets 1$ \Comment{shaking neighborhood}
\State $\text{iterCnt} \gets 1$
\Repeat
    \For{$o \in \mathcal{O}$} 
        \State $\text{prob\_operator}_o \gets \frac{\lambda_o}{\sum_{o \in \mathcal{O}}\lambda_o}$
    \EndFor
    \State Randomly select next operator $o$ based on $(\text{prob\_operator}_o)_{o\in\mathcal{O}}$
    \State $(\mathcal{R}', \phi') \gets \text{apply\_operator}((\mathcal{I},\mathcal{R}, \phi), o, T, n_{\text{max}})$
    \If{$(\mathcal{R}',\phi')\notin \mathcal{L}$ and $\text{Random}(0,1) < \exp\left(\dfrac{\text{cost}(\mathcal{R},\phi) - \text{cost}(\mathcal{R}', \phi')}{T}\right)$} 
        \State $(\mathcal{R}, \phi) \leftarrow (\mathcal{R}', \phi')$ 
        \State Update tabu list $\mathcal{L}$ respecting $N^{\text{tabu}}$ \Comment{tabu simulated annealing}
        \If{$\text{cost}(\mathcal{R},\phi) < \text{cost}(\mathcal{R}^{\text{best}}, \phi^{\text{best}})$}
            \State $(\mathcal{R}^{\text{best}}, \phi^{\text{best}}) \gets (\mathcal{R}, \phi)$
        \EndIf
        \State Select $\sigma_{\textit{iterCnt}}$ from $\sigma$ based on acceptance criterion
    \EndIf
    \State $\lambda^{\text{seg}}_o \gets \lambda^{\text{seg}}_o + \sigma_{\text{iterCnt}}$ \Comment{segment score update}
    \If{$\text{iterCnt} \text{ mod } N^{\text{seg}} = 0$} \Comment{global score update after segment end}
        \State Update global scores $\lambda$ using $\lambda^{\text{seg}}$ and $\rho$
    \EndIf
    \State $T \gets T^{\text{start}} \cdot \left(1 - \dfrac{t^{\text{elap}}} {t^{\text{max}}}\right)$ \Comment{temperature update based on elapsed runtime}
        \State $\text{iterCnt} \gets \text{iterCnt} + 1$
    \If{shaking criterion reached}
        \State $k \gets \text{UpdateNeighborhood}(k, (\mathcal{R}^{\text{best}}, \phi^{\text{best}}), (\mathcal{R}, \phi))$ \Comment{increase or reset neighborhood}
        \State $(\mathcal{R},\phi) \gets \text{shake}(\mathcal{I},(\mathcal{R}^{\text{best}},\phi^{\text{best}}),\Lambda_k, N^{\text{shaking}})$
    \EndIf
        
\Until $t^{\text{max}}$ reached
\State \Return $(\mathcal{R}^{\text{best}}, \phi^{\text{best}})$
\end{algorithmic}
\end{algorithm}
At the beginning of the algorithm, an initial feasible solution is constructed, and necessary initializations are performed. At each iteration, an operator is selected randomly based on operator scores $(\lambda_o)_{o\in\mathcal{O}}$. Note that in this context, only the facility string replacement operator presented in \Cref{sec:facility_string_replacement} utilizes $T$ and $n_{max}$.
After the operator is applied and returns a new solution $(\mathcal{R}',\phi')$, a \textit{tabu simulated annealing} (TSA) criterion is used to decide whether the new solution is accepted. After this, the current temperature is increased, and the tabu list and the operator scores are updated. To stabilize operator selection, (global) score updates are applied only after a certain number of iterations $N^{\text{seg}}$, i.e., at the end of a \textit{segment}. Whenever no feasible moves can be found, or no new best solution is found for a certain number of iterations, a shaking step is performed to move to a different region of the solution space. The algorithm terminates once the runtime limit $t^{\text{max}}$ is reached.\\
 In the following sections, we provide detailed descriptions for the individual components of the proposed algorithm. For this, we denote by $(\mathcal{R}, \phi)$ the current solution, and by $(\mathcal{R}', \phi')$ the solution obtained by applying a local search operator to $(\mathcal{R},\phi)$.

\subsection{Constructive Heuristic}\label{sec:constructive_heuristic}
The proposed constructive heuristic consists of two components. First, a feasible assignment of customers to facilities is constructed. Given an initially empty solution, a random closed facility, i.e., a facility that is not part of the current solution, is selected. All customers that can be covered by the selected facility are sorted in ascending order with respect to the number of closed facilities that can cover them. Then, the facility is opened and customers are iteratively assigned to the selected facility until its capacity is reached or all customers are assigned. This process is then repeated until all customers are assigned to exactly one facility. Second, given the selected facilities and underlying customer assignments, routes are constructed. For this, a random selected facility $i$ is picked. Then, an empty route $\{d_i, d_i\}$ is initialized, where $d_i = \text{argmin}\left\{\text{dist}_{di}: \ d\in D\right\}$ is the distance-wise closest depot to $i$. Facilities are chosen greedily based on their cheapest insertion distance, i.e., the minimum increase in route distance when inserting the facility in the current route, and then inserted at the respective cheapest insertion points, until either the vehicle capacity or the maximum facility count $p$ or route length $q$ is reached. Note that, in this step, insertion points and distances for all not-yet-assigned facilities are recomputed after every route update. This process is repeated until all facilities are part of exactly one route.\\
If each customer is assigned to a facility and each facility is assigned to a route, a feasible solution was found. On the contrary, if a customer cannot be assigned to any facility without violating its capacity, or a facility cannot be assigned to any route due to capacity or length restrictions or because all vehicles are used up, the run is considered failed. In either case, the constructive heuristic is restarted using a different random seed. This process is repeated $N^{\text{start}}\in \mathbb{N}$ times. If no feasible solution was found in any run, the AVNS terminates without a solution. Otherwise, the best found solution is returned.

\subsection{Adaptive Operator Selection}\label{sec:adaptive_operator_selection}
Because the problem encompasses several heterogeneous subproblems, such as the m-CTP and the VRDAP, different operators may perform well on different instance characteristics and at different stages of the algorithm. Therefore, an adaptive operator selection mechanism is employed to dynamically adjust the search behavior. Given an operator set $\mathcal{O}$, each operator $o\in\mathcal{O}$ is assigned a score $\lambda_o\geq 0$ that measures its past performance. If an operator frequently returns improving moves, its score tends to be higher than that of other operators. In each iteration, we randomly select an operator with a probability based on its relative score $\dfrac{\lambda_o}{\underset{o\in\mathcal{O}}{\sum}\lambda_o}$, thus biasing the selection toward favoring operators that have proven to be more efficient than others. We note that operators are not permitted to return non-changing moves; i.e., they either return no solution (when no feasible changes were found) or a solution that is distinct from the input solution.\\
Every time an operator $o\in\mathcal{O}$ returns a new solution $(\mathcal{R}',\phi')$, its \textit{segment score} $\lambda^{\text{seg}}_o \geq 0$ is updated as $\lambda_o^{\text{seg}} \gets \lambda_o^{\text{seg}} + \sigma$, where $\sigma\geq 0$ is a score increase that depends on the solution quality. Following the considerations by \cite{Sacramento2019}, we utilize a four-split scoring system $\sigma\in\left\{\sigma_{\text{global}}, \sigma_{\text{local}},\sigma_{\text{sa}},\sigma_{\text{reject}} \coloneqq 0\right\}$ as follows:
\begin{align*}
    \sigma = \begin{cases}
        \begin{aligned}
            &\sigma_{\text{global}} && \text{if} \ \  c(\mathcal{R}',\phi') < c(\mathcal{R}^{\text{best}},\phi^{\text{best}}),\\
            &\sigma_{\text{local}} && \text{if} \ \  c(\mathcal{R}^{\text{best}}, \phi^{\text{best}}) \leq  c(\mathcal{R}',\phi') < c(\mathcal{R}, \phi),\\
            &\sigma_{\text{sa}} && \text{if} \ \  c(\mathcal{R}',\phi') > c(\mathcal{R}, \phi) \ \text{and solution was accepted by TSA},\\
            &0 && \text{else.}
            \end{aligned}
    \end{cases}
\end{align*}
Score updates after every operator call can, particularly in early iterations, introduce significant instability. Furthermore, because the initial constructive solution typically has a large objective value, finding new local or global best solutions in the first few iterations  of the algorithm is comparatively easy. Hence, operators selected frequently in the first few iterations tend to receive overly optimistic scores. To mitigate these effects, we follow the approaches by \cite{Ropke2006} and \cite{Sacramento2019} and only update scores at the end of a \textit{segment}, i.e., a predefined number of iterations $N^{\text{seg}}$. At the end of a segment, we update scores as
\begin{align*}
    \lambda_o = \rho \cdot \lambda_o + (1-\rho) \cdot \dfrac{\lambda^{\text{seg}}_o}{\pi_o},
\end{align*}
where $\pi_o$ is the number of calls of operator $o$ in the current segment, and $\rho\in [0,1]$ is a smoothing factor. The smoothing factor $\rho$ allows for balancing the impact of recent and historical operator performance. Finally, we note that tabu moves do not yield a score increase, even if they would have been accepted by the simulated annealing criterion.

\subsection{Tabu Simulated Annealing}\label{sec:tabu_simulated_annealing}
After applying an operator, a new solution $(\mathcal{R}',\phi')$ is obtained. A naive approach is to accept the new solution only if it is improving, i.e., if $c(\mathcal{R}', \phi') < c(\mathcal{R}, \phi)$ holds. However, such local search-based approaches often get stuck in local minima. Accepting non-improving solutions with a certain probability can allow the algorithm to escape local minima without leaving the local search neighborhood. Therefore, following \cite{Sacramento2019}, we accept a non-improving solution with probability $\text{exp}\left(\dfrac{\text{cost}(\mathcal{R},\phi) - \text{cost}(\mathcal{R}', \phi')}{T}\right)$, where the temperature $T > 0$ is a parameter that aims to control the level of exploration and is gradually reduced over time. Furthermore, because non-improving moves can be accepted, cycles can occur. To break such cycles, we define a tabu list $\mathcal{L}$ containing the most recent accepted solutions. Every time a new solution is accepted, $\mathcal{L}$ is updated to include said solution. Whenever the length $\vert\mathcal{L}\vert$ exceeds a maximum length $N^{\text{tabu}}$, its least recent solution is removed.

\subsection{Shaking}\label{sec:shaking}
To enable exploration of areas of the solution space that are dissimilar to previously found solutions, we employ a shaking procedure. For this, let $(\mathcal{R},\phi)$ be an arbitrary solution. Shaking is performed in two distinct situations. First, when all operators return either no feasible changing moves or tabu moves since the last solution change, it is impossible to escape the current solution; hence, shaking is performed. For the second criterion, we say that $(\mathcal{R}',\phi')$ is \textit{locally improving} $(\mathcal{R},\phi)$ if $\text{cost}(\mathcal{R},\phi) > \text{cost}(\mathcal{R}',\phi')$ holds. After performing a shaking step, the shaken solution, say $(\mathcal{R}^{\text{shaking}}, \phi^{\text{shaking}})$, is stored. Whenever a solution $(\mathcal{R}', \phi')$ that is locally improving $(\mathcal{R}^{shaking}, \phi^{shaking})$ is found, we update
\begin{align*}
    (\mathcal{R}^{\text{shaking}}, \phi^{\text{shaking}}) \gets (\mathcal{R}',\phi').
\end{align*}
If $(\mathcal{R}^{\text{shaking}}, \phi^{\text{shaking}})$ is not updated in the last $N^{\text{stag}}$ iterations, that is, no new local best solution was found in the last $N^{\text{stag}}$ local search operator calls, the solution space around the current solution is considered sufficiently explored, triggering a shaking step to explore other regions of the solution space.\\
The shaking procedure proposed in this paper follows a destroy-and-repair paradigm. In each shaking step, percentages of all customers assignments, selected facilities and depots that are to be removed are defined. These percentages thus control the degree of exploration and are gradually increased.\\
In the first step, parts of the current solution are removed to obtain a skeleton solution. In particular, let $(d=v_0,v_1,\dots,v_l,v_{l+1}=d)$ be a route in the current solution with $l \in\mathbb{N}$ and $d\in D$. Then, for each $i \in \{1,\dots,l\}$, we calculate a \textit{removal delta} as 
\begin{align*}
    \Delta_{v_i} \coloneqq c_{v_{i-1}, v_{i+1}} - c_{v_{i-1}, v_i} - c_{v_i, v_{i+1}} + \underset{l\in W: \phi(l) = v_i}{\sum}\left(\underset{f \in F_l}{\min}\left\{a_{l, f}\right\} - a_{l,v_i}\right).
\end{align*}
Note that we assume $c_{dd} = 0$. The removal delta of a facility $v_i$ aims to estimate the (maximum) potential saving obtained when removing $v_i\in\mathcal{F}$. Facilities are iteratively randomly selected for removal based on a probability-defining softmax function $\sigma: \mathcal{F}\rightarrow [0,1]$, until the desired percentage of facilities is removed:
\begin{equation*}
    \sigma(i) \coloneqq \dfrac{\exp\left({-\dfrac{\Delta_{v_i}}{\Delta_{\text{max}}}}\right)}{\underset{j\in\mathcal{F}}{\sum} \exp\left({-\dfrac{\Delta_{v_j}}{\Delta_{\text{max}}}}\right)},
\end{equation*}
where $\Delta_{\max} \coloneqq \max\left\{\vert \Delta_{v_i} \vert: \ v_i\in\mathcal{F}\right\}$ holds. This way, the procedure is biased toward favoring the removal of potentially inefficient facilities, while still maintaining the intrinsic random character of a shaking procedure. Note that the normalization using $\Delta_{max}$ puts additional emphasis on exploration. Furthermore, note that, whenever a facility is removed, all of its customers become unassigned. If, after removing the required number of facilities, additional customer coverings need to be removed, an analogous procedure is applied to remove additional customer coverings using removal deltas $\Delta_{l} \coloneqq \underset{f\in F_l}{\min} \{a_{l,f}\} - a_{l,\phi(l)}$. The same holds for depots, whereas removal deltas are calculated as
\begin{align*}
    \Delta_d \coloneqq \underset{R \in \mathcal{R}(d)}{\sum} \ \underset{v\in R}{\sum} \  \underset{l\in W: \phi(l) = v}{\sum} r_{l}
\end{align*}
where $\mathcal{R}(d)$ is the set of routes that start and end at depot $d\in D$, and $Q(R)$ is the total demand served by route $R$.\\
Following the removal of customer assignments, facilities, and depots, an attempt is made to extend this infeasible skeleton solution into a complete feasible solution using analogous techniques as described in \Cref{sec:constructive_heuristic}. If the extension is successful, the shaking procedure terminates by returning the perturbed solution. Otherwise, the process is repeated up to $N^{\text{shaking}}$ times. If no feasible solution is obtained within the prescribed number of runs, the procedure outputs the original input solution.\\
Before each shaking call, the shaking neighborhood is updated. If the best found solution $(\mathcal{R}^{\text{best}}, \phi^{\text{best}})$ was updated since the last shaking call, the neighborhood is reset to $k=1$. Otherwise, it is increased by one, i.e. $k\gets k+1$. If $k=K$ holds, i.e., the last neighborhood was used in the previous shaking step, it is reset to $k=1$.

\subsection{Local Search} \label{sec:local_search}
We employ several types of operators to account for the fact that the problem at hand incorporates routing and covering decisions at once. To improve route structures, we rely on multiple well-established intra-route (2opt, swap, relocate) and inter-route (1-point move, 2-point move, 2-opt*, cross-string) operators known from the literature. For a detailed description of these operators, we refer to \cite{groer2010library}, \cite{funke2005local}, and \cite{allahyari2015hybrid}. The implementation of said operators is adapted from \cite{Rasku2019}. All of the aforementioned operators are searched in a best-apply fashion, that is, all possible moves of an operator are evaluated, and the best move is returned.\\
In the following, we propose a greedy assignment heuristic to improve assignments of depots to routes, as well as an operator that iteratively removes facility strings and reassigns customers to already opened or currently closed facilities. \\
We note that local search operators typically return the input solution if no improving move is found. However, because non-improving moves are permitted in the proposed solution approach, all operators are implemented such that they never return a solution that is equal to the input solution.

\subsubsection{Greedy Depot Assignment}\label{sec:depotIP}
For an optimal reassignment of routes to depots, we make use of a greedy assignment heuristic. For each route $R \coloneqq (d=v_0,v_1,\dots,v_l,v_{l+1}=d)\in \mathcal{R}$ and each depot $d^*\in D$, we then compute the cost delta obtained when first removing $d$ from $R$ and then inserting $d^*$ at its cheapest position.
We then select $\bar{d}_{R} \coloneqq \text{argmin}\left\{\Delta c_{R,d^*}: d^*\in D\right\}$, reassign $R$ to $\bar{d}_{R}$ and insert $\bar{d}_{R}$ at its cheapest position. To ensure that this operator does not return the input solution, we check if said reassignment changed the depot assignment of at least one route. If that is not the case, we compute the second-best depot $\underline{d}_{R} = \text{argmin}\left\{\Delta c_{R,d^*}: \ d^*\in D \setminus\{\bar{d_R}\}\right\}$ for each route $R$. Let $\underline{R}$ be the route that minimizes $\left\vert \Delta c_{\underline{R},\underline{d}_{\underline{R}}} -\Delta c_{\underline{R},\bar{d}_{\underline{R}}}\right\vert$. We then alter the assignment by assigning $\underline{R}$ to $\underline{d}_{\underline{R}}$. All other assignments remain unchanged.

\subsubsection{Customer Assignment: Facility String Replacement}\label{sec:facility_string_replacement}
In the following, we propose an operator that iteratively removes strings of facilities from routes and then reassigns their customers to open or closed facilities, thus simultaneously modifying route structure and customer assignments. Unlike other operators that focus on routing decisions, this operator explicitly accounts for the interaction between facility selection and assignment decisions. Its pseudocode is presented in \Cref{algo:GreedyStringReplacement}.
\begin{algorithm}[H]
\caption{FacilityStringReplacement$(sl, \mathcal{R}, \phi, \mathcal{F}, n_{\text{max}}, T)$}\label{algo:GreedyStringReplacement}
\begin{algorithmic}[1]
\Require String length $sl\in\mathbb{N}_{>0}$, current solution $(\mathcal{R}, \phi)$, open facilities $\mathcal{F}$, maximum number of accepted solution changes $n_{\text{max}}$ per facility, current temperature $T$
\State $U \leftarrow \mathcal{F}$ \LineComment{list of unchecked facilities}
\State $n_f \leftarrow 0 \ \forall f \in F$ \LineComment{accepted move counter per facility}
\While{$U \neq \emptyset$}
    \State $f \leftarrow \text{random facility in } U$; $\ R \leftarrow$ route containing $f$; $\ i \leftarrow$ index of $f$ in $R\in\mathcal{R}$
    \State $U \gets U\setminus\left\{f\right\}$
    \If{$n_f \geq n_{\max}$}
        \State \textbf{Go to} Step 4 \LineComment{$f$ has been the start of too many accepted solution changes}
    \EndIf
    \State $S \leftarrow \{v_i, \ldots, v_{\min\{\vert R \vert - 2,\, i+sl-1\}}\}$
    \State $\mathcal{W} \leftarrow \phi^{-1}(S)$ \LineComment{customers to reassign}
    \State $S_{\text{new}}\leftarrow \emptyset$; $\ \phi_{\text{new}} \leftarrow \phi$ \LineComment{new string and customer assignments}
    \State Re-insert all mandatory facilities from $S$ into $S_{\text{new}}$ according to their order in $S$
    \State Compute scores $sc(l, f')$ for all $l \in \mathcal{W}$, $f' \in F_l$
    \State $\mathcal{H} \gets \left\{(l, f'): \ l \in \mathcal{W}, f'\in F_l\right\}$ \Comment{all possible reassignments}
    \While{$\mathcal{H} \neq \emptyset$ \textbf{and} $\mathcal{W} \neq \emptyset$}
        \State $(l, f') \leftarrow \text{argmin}\left\{sc(l,f'): \ (l, f') \in \mathcal{H}\right\}$ \Comment{pick most promising reassignment first}
        \State $\mathcal{H} \gets \mathcal{H} \setminus \left\{(l, f')\right\}$
        \If{assigning $l$ to $f'$ is feasible w.r.t.\ $Q$, $Q^F$, $p$, and $q$}
            \State $\phi_{\text{new}}(l) \leftarrow f'$; $\ \mathcal{W} \leftarrow \mathcal{W} \setminus \{l\}$
            \State Remove $(l, f^*)$ from $\mathcal{H}$ for all $f^*\in F_l$
            \If{$f' \notin (\mathcal{F}\setminus S)\cup S_{\text{new}}$}
                \State Open $f'$, insert into $S_{\text{new}}$ at cheapest position \LineComment{new facility opened}
                \State Update scores $sc(l^*, f')$ for all $l^*\in \mathcal{W}$ with $(l^*,f')\in\mathcal{H}$
            \EndIf
        \EndIf
    \EndWhile
    \If{$\mathcal{W} \neq \emptyset$}
        \State \textbf{Go to} Step 4 \LineComment{not all customers could be reassigned}
    \EndIf
    \State Evaluate solution acceptance using simulated annealing with temperature $T$
    \If{solution is accepted}
        \State Replace $S$ with $S_{\text{new}}$ in $R$; update $\mathcal{R}$ and $\mathcal{F}$ accordingly
        \State $\phi \leftarrow \phi_{\text{new}}$
        \State $n_f \leftarrow n_f + 1$
        \State $U\gets U \cup S_{\text{new}}$
    \EndIf
\EndWhile
\State \textbf{Return} new solution $(\mathcal{R}, \phi)$
\end{algorithmic}
\end{algorithm}

The operator begins by initializing a list $U$ consisting of all open facilities for which subsequent strings of length $sl$ are to be attempted to be removed. Then, a random facility $f\in U$ and the corresponding string $S$ of length $sl$, whose first facility is $f$, is selected and the corresponding facilities are removed from its route $R=(d,v_1,\dots,v_{\vert R \vert - 2}, d\}$. By doing so, all customers that were assigned to a facility in $S$ become unassigned. A new string $S_{\text{new}}$ is initialized, which is subsequently populated with newly opened facilities. Mandatory facilities, i.e., facilities that cover a customer that can not be covered by any other facility or facilities that must be visited, are immediately re-inserted into $S_{\text{new}}$ in the same order as in $S$. Then, given an unassigned customer $l\in \phi^{-1}(S)$ and a facility $f'\in F_l$, we compute a score $sc(l,f')$ as
\begin{equation*}
    sc(l, f') \coloneqq \begin{cases}
        a_{l,f'} - a_{l,\phi(l)} & \text{if } f'\in (\mathcal{F} \setminus S)\cup S_{\text{new}},\\
        \dfrac{r_l}{Q^F_{f'}}\cdot \Delta^R_{f'} + a_{l,f'} - a_{l,\phi(l)} & \text{else},
    \end{cases}
\end{equation*}
where $\Delta^R_{f'}$ is the routing cost change obtained from inserting $f'$ at its cheapest position into $S_{\text{new}}$, taking into account that $S_{\text{new}}$ is to be placed at the same position where $S$ was previously removed. If a facility $f'$ is already open in the current solution, the score $sc(l, f')$ equals the covering cost changes when assigning $l$ to $f'$. If a facility $f'$ is currently closed, the score also considers a share of the change of routing costs incurred when inserting $f'$ into the current new string $S_{\text{new}}$ at its best position.\\
All pairs $(l,f')$, which have not been evaluated for reassignment yet, are then stored in a set $\mathcal{H}$. Iteratively, the tuple with the lowest score is selected from $\mathcal{H}$, and it is verified if assigning the respective customers to the respective facility is feasible with respect to capacity restrictions, as well as to the maximum facility count and route length constraints (\ref{MIP:p-constr2})--(\ref{MIP:q-constr3}). Note that the latter two conditions are only checked when the target facility is currently closed. If the assignment is feasible, all relevant data structures, i.e., $S_{\text{new}}$, $\phi_{\text{new}}$, $sc(l^*,f')$, $\mathcal{W}$, and $\mathcal{H}$ are updated. The process terminates when all tuples in $\mathcal{H}$ have been evaluated, or all customers have been assigned, i.e. when $\mathcal{H} = \emptyset$ or $\mathcal{W} = \emptyset$ hold. If $\mathcal{W} \neq \emptyset$ holds, the string replacement is considered failed, all changes are reverted, and the next string is checked. Otherwise, we utilize the proposed simulated annealing criterion to assess whether the changed solution should be accepted. If the solution is accepted, $U$ is updated to include all newly added facilities, $S_{\text{new}}$ is introduced into $R$ at the prior position of the removed string $S$, $\phi$ is updated based on $\phi_{\text{new}}$, and the next string is selected.
The operator terminates once $U$ is empty or once every facility was selected $n_{\text{max}}$ times as the starting point of a to-be-evaluated string. Preliminary computational results showed that setting $n_{\text{max}} = 3$ provides an optimal trade-off between performance and computational effort. Note that we call FacilityStringReplacement with a string length $sl\in\{1,3,5\}$. In particular, for each $sl\in\{1,3,5\}$, we include a separate operator $\text{FacilityStringReplacement}(sl, \mathcal{R}, \phi, \mathcal{F}, n_{\text{max}}, T)\in \mathcal{O}$ into the set of available operators.

%% file: sections/computational_experiments.tex
\section{Computational Experiments}
\label{sec:comp_experiments}
The computational study is designed to assess the effectiveness of the proposed AVNS along two dimensions: (i) its ability to produce high-quality solutions across different SCRP variants, and (ii) its robustness when applied to problem instances with distinct structural characteristics. In particular, we investigate whether the proposed AVNS can achieve competitive performance compared to specialized state-of-the-art algorithms for the VRDAP \citep{Oliveira2025_VRDAP_ILS} and the m-CTP \citep{Oliveira2025_mCTP_exact}, while maintaining general applicability across problem classes.
To this end, we apply the proposed AVNS to benchmark instances from the VRDAP and the m-CTP, using the same instance sets as in the aforementioned studies to ensure a fair comparison.\\
The proposed algorithm is implemented in Python 3.12. All computational experiments are conducted on a standard workstation equipped with an Intel Core i7-8700K 12-core 3.7GHz processor, 32GB of RAM and running Windows 11. To ensure reproducibility, all benchmark instances, solution data, and parameter configurations can be found under \url{https://github.com/andreashagntum/RoutingCovering_Instances}. The implementation of the proposed AVNS is accessible under \url{https://github.com/andreashagntum/RoutingCovering_AVNS}.\\
In the remainder of this section, we first describe the parameter tuning procedure, followed by the presentation of computational results.

\subsection{Parameter Tuning\label{sec:parameter_tuning}}
To calibrate the parameters of the proposed AVNS, we employ a Bayesian optimization approach using the Python package SMAC3 (cf. \cite{lindauer-jmlr22a}). For this purpose, we select 10 representative test instances from each instance set. We then parallelize the framework implemented in the aforementioned package to 4 threads and apply it for a maximum search duration of 24 hours, solving each instance 5 times per trial, i.e., per evaluated parameter configuration. This results in a total of 41 trials for each m-CTP instance set, and 277 for the VRDAP instance set. 
The performance of each parameter configuration was measured by computing the optimality gap of the best solution found and averaging them over all tuning instances. For all instance classes except VRDAP and case study instances, lower bounds reported in the literature were used. For the VRDAP and the Afghanistan case study instances, we computed the lower bounds ourselves using the LP relaxation of the MILP formulation presented in \Cref{subsec:MIP}. To obtain a tighter relaxation, we strengthened the model by incorporating the following additional constraints: 
\begin{align*}
      &\sum_{(j, i) \in A}f_{ji} - \sum_{(i, j) \in A}f_{ij} = w_i \quad \text{for all } i \in F, \\ 
    &0 \leq f_{ij} \leq M x_{ij} \quad \text{ for all } (i, j) \in A,  
\end{align*}
with new continuous variables $f_{ij} \text{ for } (i, j) \in A$, whereas the parameter $M$ can be set to $p$ if $p < \infty$ or to $|F|$ otherwise.\\
Furthermore, we adjust the maximum runtime $t_{\max}$ (measured as CPU runtime) based on the runtime of the state-of-the-art algorithms for the considered problem sets. For all m-CTP(-p) problem variants, \cite{Oliveira2025_mCTP_exact} provides the best-performing algorithm. It utilizes a heuristic with a runtime of 180 seconds to obtain initial upper bounds, which are then used in their branch-and-price approach with a runtime limit of 4800 to 7200 seconds, depending on the instance sets. Because the aforementioned authors use hardware that is comparable to ours, we thus adapt a runtime limit of $t_{\max} = 180$ seconds per run. For the VRDAP instances, \cite{Oliveira2025_VRDAP_ILS} can be considered the state of the art, applying an iterative local search up to 100 times with 100 iterations per call, resulting in an average runtime of 18 to 33 seconds, depending on the instance sizes. According to \cite{cpu_benchmark}, the CPU used in the following achieves performance results on single-thread benchmark tests that are 18\% better than those of the hardware used by \cite{Oliveira2025_VRDAP_ILS}. We thus impose a runtime limit $t_{\max} = 25$ seconds per run.\\
Moreover, in all conducted experiments, we fix the following parameters, because they have shown to only have a minor impact on performance: $N^{\text{start}} = 200$, $N^{\text{shaking}} = 20$, and $\lambda^{\text{initial}} = (10)_{o \in \mathcal{O}}$. Finally, for each instance set and the subsequent computational experiments, we select the best configuration found by SMAC.

\subsection{Results on Benchmark Instances}\label{sec:performance_comparison}
We collected instances originally proposed by \cite{Ghoniem2013} (VRDAP), \cite{Pham2017} (m-CTP), \cite{Glize2020} (m-CTP-p), and two instance sets by \cite{Oliveira2025_mCTP_exact} (m-CTP-p). We applied our algorithm to a total of 300 VRDAP test instances, 192 m-CTP instances, and 212 m-CTP-p instances. Note that 58 of 212 m-CTP-p instances have mandatory facilities. In the following, we denote the test instance sets as \textit{VRDAP (Ghoniem)}, \textit{m-CTP (Pham)}, \textit{m-CTP-p (Glize)}, \textit{m-CTP-p (Oliveira)}, and \textit{M-m-CTP-p (Oliveira)}, where the prefix $M$- indicates that the test instances contain mandatory facilities. The names in brackets indicate the authors of the papers in which the respective instance set was first proposed. Note that in the aforementioned m-CTP(-p) instances, the parameter $q$ limits the maximum routing cost of each route. Each instance is solved 5 times. In the following, we report performance for the best, median, and worst solutions found for each instance.\\

\begin{table}[h]
\centering
\begin{tabular}{llcc}
\toprule
Instance Set & \makecell[l]{Inst. with \\ no / feas. / opt.\ sol.} & \makecell[c]{Avg.\ UB gap\ \\ (best / median / worst sol.)} & \makecell[c]{\# new \\ feas.\ / best sol.} \\
\midrule
m-CTP (Pham) & 0 / 1 / 191 & 0.13\% / 0.49\% / 0.80\% & 1 / 0 \\
m-CTP-p (Glize) & 0 / 0 / 96 & 0.07\% / 0.19\% / 0.28\% & 0 / 0 \\
VRDAP (Ghoniem) & 0 / 300 / 0 & 0.59\% / 0.98\% / 1.44\% & 0 / 34 \\
m-CTP-p (Oliveira) & 30 / 7 / 21 & 2.92\% / 4.74\% / 6.25\% & 30 / 0 \\
M-m-CTP-p (Oliveira) & 35 / 5 / 18 & 1.38\% / 2.09\% / 3.33\% & 35 / 0 \\
\bottomrule
\end{tabular}
\caption{Summary of computational results by instance set.}
\label{tab:results_summary}
\end{table}

\Cref{tab:results_summary} summarizes the performance of the proposed algorithm on all benchmark sets. Column \textit{Inst. with no / feas / opt.sol.} denotes the number of instances, for which no feasible solution is known, a feasible but not necessarily optimal solution is known, and the optimal solution is known in the literature, respectively. Column \textit{Avg. UB gap} shows the gap, computed as $\dfrac{UB_{\text{heur}} - UB_{\text{lit}}}{UB_{\text{lit}}}$, where $UB_{\text{heur}}$ is the objective value of the best (resp. median, worst) solution found by the proposed heuristic across all 5 runs, and $UB_{\text{lit}}$ is the value of the best known solution from the literature. For this, only instances for which an upper bound exists in the literature are considered. Moreover, column \textit{\# new feas. / best sol.} contains the number of instances for which we are the first to propose a feasible solution, and the number of instances for which we propose a new best solution, respectively. Note that the latter value does not include instances for which no feasible solution was previously known. For the instance sets m-CTP (Pham) and m-CTP-p (Glize), the proposed algorithm achieves stable and near-optimal performance, considering the upper bounds from the literature already comprise the optimal solutions for all said instances. These results are particularly noteworthy given that our algorithm operates within a runtime limit of 180 s, whereas the state-of-the-art approach of \cite{Oliveira2025_mCTP_exact} has a runtime limit of 4800 s and average runtime of 240 s (plus 180 s for the initial upper bound computation) per instance.
For the VRDAP instances, gaps between the best and worst run range between 0.59\% and 1.44\%. Notably, the algorithm finds new best solutions for 34 out of 300 VRDAP instances, despite operating under a runtime limit of only 25 seconds. Recall that said limit is comparable to the average runtime of the iterative local search of \cite{Oliveira2025_VRDAP_ILS}, which is specifically tailored to this problem class. 
For the instance sets (M-)m-CTP-p (Oliveira), the computed gaps increase from 1.38\% and 2.92\% to 3.33\% and 6.25\%, respectively, indicating a slight performance and robustness degradation. This is primarily attributable to the significantly larger size of these instances relative to the other benchmark sets, as they contain up to 392 customers, 392 facilities, and 94 vehicles, while all other instance sets contain at most 149 customers, 99 facilities, and 11 vehicles. Importantly, however, the proposed heuristic finds feasible solutions for all instances across all sets, including all 66 instances for which no feasible solution was previously known in the literature.
Thus, the proposed AVNS remains competitive with problem-specific state-of-the-art algorithms on all instance sets in a fraction of the runtime of the benchmark approaches, highlighting its robustness and its ability to handle structurally different problem variants. In addition, it finds feasible solutions for every test instance, demonstrating the reliability of the algorithm.

%% file: sections/case_study.tex
\section{Case Study: Disaster Relief Operations in Afghanistan}
\label{sec:case_study}
\begin{figure}[ht]
\centering
\includegraphics[scale=0.4]{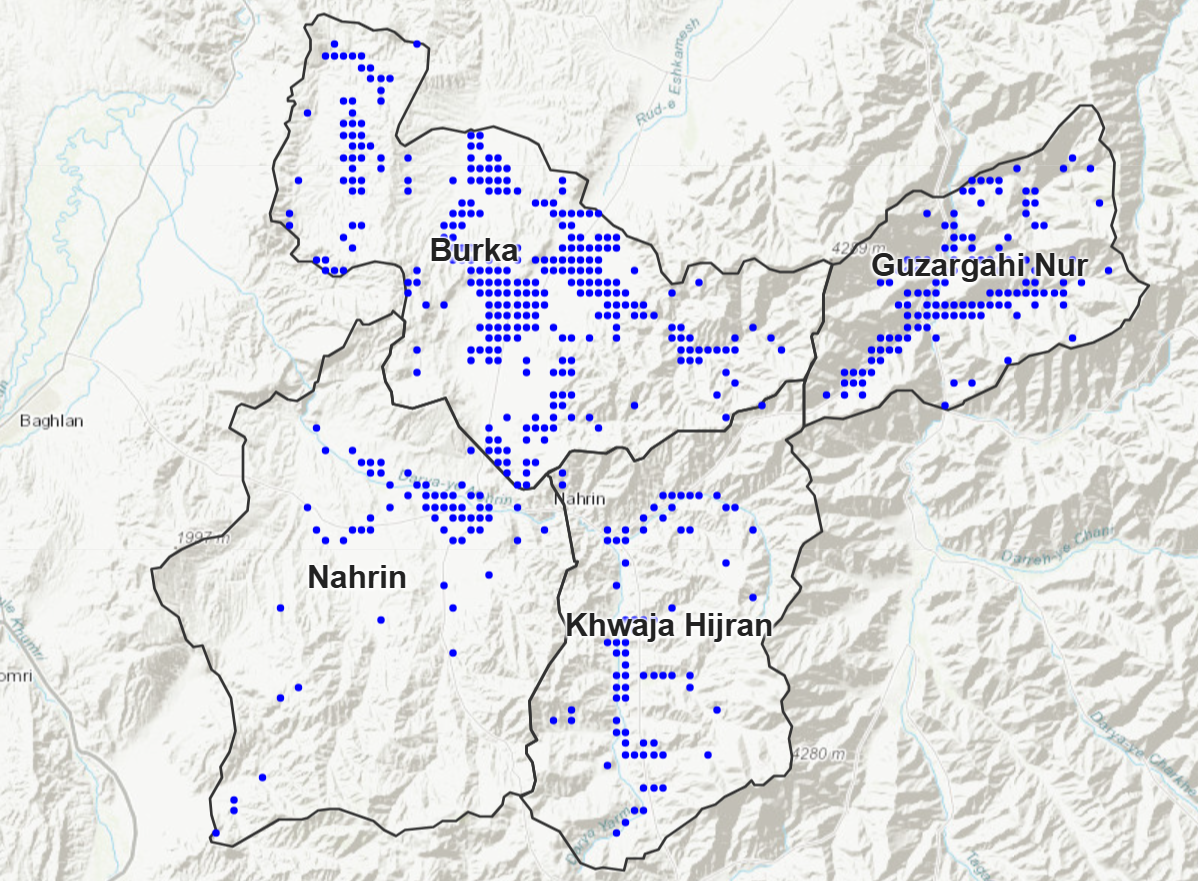}
\caption{Considered affected region}
\label{fig:instance_setting}
\end{figure}

In the following, we apply the algorithm proposed in \Cref{sec:algorithm} to derive a post-acute food aid response strategy for the May 2024 flash floods in Afghanistan. We focus on a total of 4 districts in the province of Baghlan, which were particularly heavily affected and inaccessible from the outside for several days as access routes were flooded or destroyed, while transport by vehicles within the aforementioned regions was also largely impossible (cf. \cite{inaccessible_source_who}, \cite{inaccessible_source_unicef}, \cite{awd_cholera}). \Cref{fig:instance_setting} shows the considered region. Each blue dot represents an area of size 1 km², in which at least one individual that was affected by the May 2024 flash floods is located. The 4 considered districts are separated by black boundaries and labeled accordingly.\\
Because a significant portion of the Afghan population is malnourished, a consistent supply with nutrients is crucial to reduce mortality rates. In humanitarian aid distribution to each individual's location is often not practically feasible (cf. \cite{RANCOURT201568}, \cite{NajiAzimi2012}). Hence, we identify helicopter landing sites that are expected to be usable even in adverse conditions. We assume that of each affected household, one person travels on foot to a designated landing site to collect their aid. Moreover, we integrate routing distances as a measure for financial expenses, resulting from chartering, operation and upkeep of helicopters, as well as burden put on the affected individuals (in the following also called beneficiaries) resulting from ground travel efforts, into one objective function.\\
Several factors, such as the number of beneficiaries, demands, fleet size, and feasible beneficiary-to-landing site assignments, influence the effectiveness of a potential response strategy. The purpose of the following study is to derive guidelines for practitioners on when the proposed distributional structure is suitable, if and when the chartering of additional helicopters makes sense, how operational costs and burdens put on the beneficiaries by the proposed covering approach trade off, and to which amount neglecting the supply of aid to beneficiaries that are particularly hard to reach can reduce the operational costs of the response. By that, the proposed approach enables humanitarian aid organizations to devise suitable, efficient, and fair distributional structures in post-acute flood emergency situations as observed in Baghlan in May 2024.\\
In \Cref{sec:case_study_instance_generation}, we present the process of test instance construction in a condensed manner. For an in-depth explanation of the construction of test instances, including specific threshold values, we refer to the supplementary material of the paper. \Cref{sec:instance_characteristics} summarizes several key characteristics of the considered instance set and gives first insights on the limitations of the proposed distributional structure. In \Cref{sec:additional_heli_impact}, the impact of increasing fleet sizes is analyzed and recommendations for optimal maximum numbers of helicopters are presented. \Cref{sec:tradeoff_routing_vs_covering} focuses on the trade-off between operational cost and burden caused by indirect supply of beneficiaries through the aforementioned covering approach. Finally, \Cref{sec:triage} elaborates on the impact of an introduction of a service threshold, i.e., neglecting beneficaries that are hard to reach, and provides insights on when such an approach can be of interest for practitioners. All instances are publicly available under \url{https://github.com/andreashagntum/RoutingCovering_Instances}.

\subsection{Instance Generation}\label{sec:case_study_instance_generation}
Each instance $i$ is uniquely characterized by the factor of beneficiaries $N_i^{\text{ben}} > 0$, the maximum covering duration $CD_i^{\max} > 0$, the service threshold $ST^{\max}_i > 0$, the maximum number of helicopters $M_i \in \mathbb{N}$ (in the following also called \textit{fleet size}), and a routing cost scaling factor $\rho_i \geq 0$. In the following, we elaborate on the meaning of said parameters and their usage during the instance construction process.\\
We construct a given test instance $i$ as follows. First, we generate locations of affected individuals based on population densities and historical data on the number of individuals affected by the May 2024 flash floods per district, obtained from \cite{population_data} and \cite{affected_individuals}. The former reports a total of 5831 affected individuals. Then, we compute a demand (in kg of aid) for each beneficiary, i.e., for each affected individual. For this, the minimum daily calorie intake (cf. \cite{min_calorie_intake}) for each beneficiary is assumed to be covered by high-energy biscuits (HEBs), as the UN's immediate response to the considered flood event also relied heavily on HEBs (cf. \cite{wfp_heb_response}), which translates to 467 g of HEBs per beneficiary. To reduce the instance complexity, we aggregate all beneficiaries located in a square of size 1 km² into a single customer node, which inherits the total demand for all aggregated beneficiaries. To increase the diversity of the instance set and avoid conclusions that are contingent on a single population estimate, the beneficiary counts reported in \cite{affected_individuals} are then multiplied by an instance-dependent factor $N^{\text{ben}}_i$, yielding instances of varying demand intensity around the reported baseline.\\
Next, we identify landing sites, i.e., facilities, within and surrounding the affected districts. Unsuitable landing sites are omitted by analyzing elevation and terrain data by the DLR (cf. \cite{tandem_elevation_data}) and the European Space Agency (cf. \cite{esa_terrain_data}). Using data on Afghan road networks by the \cite{road_networks}, we ensure that each landing site is sufficiently close to a nearby road with no impassable obstacles, such as water bodies or mountains, in between. Because the resulting set of candidate locations is extensive and many locations are near identical, we then select a subset of 500 landing sites as facilities using a $k$-means clustering algorithm. Because many landing sites are often in close proximity to each other and do not fundamentally differ with respect to the beneficiaries they can cover, we then remove all landing sites for which another landing site is within a 12-minute walking radius.
Each beneficiary is assumed to be coverable by a landing site if road travel by foot between the landing site and beneficiary does not take longer than $CD^{\max}_i$ hours.\\
Ground travel is assumed to be undertaken by foot at a speed of $4$ km/h. Because roads might be partially unusable and require circumvention, we increase the aforementioned road distances by a factor of $1.2$. Note that this implies that covering durations can be transformed into covering distances (in km) by multiplying the covering duration by 4 and dividing the result by 1.2. All beneficiaries whose travel duration to their closest landing site exceeds $ST^{\max}_i$ hours are removed from the set of beneficiaries. At each landing site, a group of local volunteers receives the delivered goods and distributes it among the beneficiaries. Here, we assume that for each household, one beneficiary travels to their assigned landing site and collects their aid. In accordance with \cite{afg_household_data}, a household is assumed to consist of 6 individuals. To limit the work required by said volunteers and to avoid single points of failure in the distributional structure, we limit the capacity of each landing site to 200 kg.\\
Moreover, we identify all major, regional, or local airports and heliports that do not lie in the affected districts and have a distance of at most 150 km to at least one beneficiary as depot candidates. Due to insufficient data accuracy regarding the cost structures associated with airport usage in Afghanistan, particularly for humanitarian air services such as the \textit{United Nations Humanitarian Air Service} (UNHAS), we adopt a simplifying assumption and exclude depot-related cost and capacity considerations from this study. Therefore, we can also remove all depots which are not the closest depot to any landing site. The remaining 2 depots, Taloqan Airport and Bagram Airfield, are then used as potential starting points for helicopter routes.\\
A fleet of $M_i$ Mil Mi-8 AMT helicopters is available for use. Such helicopters have previously been operated by the UNHAS in Afghanistan (cf. \cite{WFP2024UNHAS}). This helicopter type is particularly suited for missions in rough circumstances, which makes it a good fit for operations in the mountainous, high-elevation areas of the considered districts. Technical details, such as range, capacities, and maximum flight height, were collected from \cite{heli_datasheet}. Additionally, we set the maximum number of stops per route to $p = 28$ in order to avoid excessively long routes. This stop count is derived by assuming $10$ minutes per stop for descent, landing, unloading, departure, and ascent, an average travel speed of $170$ km/h, a helicopter range of $q = 570$ km, and a maximum operating time of $8$ hours per helicopter.\\
Routing and covering costs are defined based on helicopter and household-weighted travel distances, respectively, and are rescaled to account for several properties. First, helicopter chartering and upkeep, as well as crews, are particularly expensive. Thus, in our modeling the cost of an additional helicopter always exceeds the potential savings realized by splitting an existing route into two routes (and thus, reducing routing distances). Second, routing costs (in the sense of objective function costs) between landing sites are proportional to the distance traveled. Third, covering costs are equal to household-weighted ground distances (and by that, directly proportional to ground travel times and the number of beneficiaries), which act as a proxy for burdens put on the beneficiaries by having to travel across rough terrain to collect their aid. Finally, routing costs are rescaled by $\rho_i$. Note that, by increasing $\rho_i$, efficient routing, and thus a reduction of financial expenses is incentivized, while reducing $\rho_i$ puts additional emphasis on short ground travel durations, thereby impacting ground response times and reducing efforts undertaken by beneficiaries.\\
Note that, as described at the beginning of this section, there is a discrepancy between customer nodes and beneficiaries. While the proposed instances were solved using said aggregation into customer nodes, all relevant performance indicators evaluated in the following are based on the number of affected households. Furthermore, unless otherwise stated, all values reported in the following are obtained by averaging the respective value over the relevant set of instances for which a solution was obtained. Finally, in the following we compute several KPIs based on covering distances rather than covering durations. As mentioned before, both values can be converted using the ground travel speed of $4$ km/h and the road distance factor of $1.2$.

\subsection{Instance Characteristics}\label{sec:instance_characteristics}
We construct test instances by varying $N^{\text{ben}}_i\in \{0.5, 0.75, 1, 1.25, 1.5\}$, $ST^{\max}_i \in \{1, 2, 3, 4\}$, $CD^{\max}_{i}\in \{1, 2, 3, 4\}$, $M_i\in\{1, 2, 3, 4\}$, and $\rho_i\in \{0.1, 1, 3, 5, 7\}$. For this, we first construct a baseline instance as elaborated on in \Cref{sec:case_study_instance_generation} by setting $M_i = 4$, $ST^{\max}_i = 4$, $CD^{\max}_i = 4$, $N^{\text{ben}}_i = 1$, and $\rho_i = 1$. For each of the aforementioned parameter combinations, we read the baseline instance and adjust it as needed. This results in 276 to 342 landing sites per instance. Note that, for all affected individuals, there exists at least one landing site within a maximum travel duration of 4 hours from their location; thus, no individuals had to be excluded from the set of beneficiaries in the baseline instance. Furthermore, instances with $ST^{\max}_i > CD^{\max}_i$ are discarded because, in such instances, at least one beneficiary is uncoverable. Of the remaining 1,000 test instances, 10 were randomly selected to perform a hyperparameter tuning as elaborated in \Cref{sec:parameter_tuning} with a runtime of 48 hours, resulting in a total of 49 tried parameter configurations. Then, all instances were solved using the AVNS proposed in \Cref{sec:algorithm} with 5 runs per instances, each with a runtime of 5 minutes. Afterward, the best solution (w.r.t. objective values) is kept for each instance and used for the following evaluations.\\
The total number of customer nodes ranges between 432 and 558, depending on $ST^{\max}_i$, which translates to 5369 to 5831 beneficiaries. After scaling with $N^{\text{ben}}_i$, the total number of modeled beneficiaries is between 2819 and 8922.\\
\begin{figure}[ht]
\centering
\includegraphics[scale=0.4]{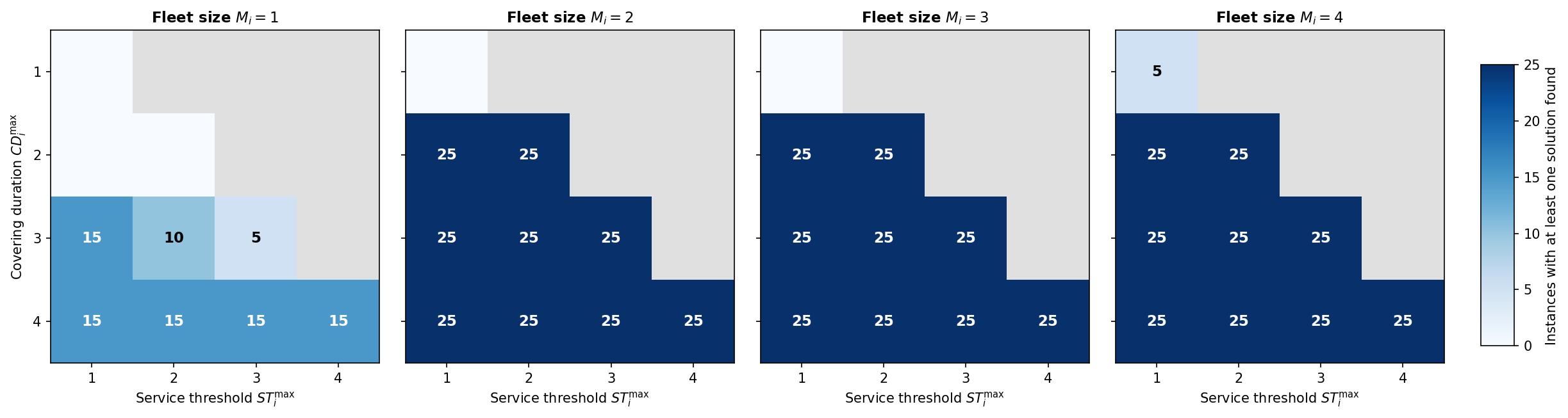} 
\caption{Heat map of instances for which a solution was found}
\label{fig:heatmap_feas_inst_cnt}
\end{figure}
\Cref{fig:heatmap_feas_inst_cnt} shows a heatmap of the number of instances for which a solution was found, filtered by all possible combinations of $M_i$, $ST^{\max}_i$, and $CD^{\max}_i$. Cells corresponding to parameter combinations that are guaranteed to be infeasible are highlighted in grey. Cells for which no feasible solution has been found for any of the underlying instances are highlighted in white. Note that each cell represents 25 instances.\\
It can be observed that when only a single helicopter is available, no solution was found for any instance with a maximum covering duration of 2 hours or less. 
Furthermore, for $CD_i^{\max}=3$ and $ST_i^{\max}=3$, the algorithm returned a solution in 5 out of 25 instances. To be precise, only instances for which $N_i^{\text{ben}} = 0.5$ holds returned a solution. For instances with $CD_i^{\max}=3$ and $ST_i^{\max} = 2$ or $ST_i^{\max}=1$, solutions were also obtained for instances with $N_i^{\text{ben}}=0.75$ and $N_i^{\text{ben}}=1$, yielding a total of 10 and 15 feasible instances, respectively. Hence, one can conclude that one helicopter is only sufficient for emergency response in medium demand scenarios, i.e., when $N^{\text{ben}}_i \leq 1$, and if covering (i.e. ground travel) durations of 3 or more hours are acceptable. At the same time, a fleet of 2 to 3 helicopters permits the use of covering durations of 2 hours, while 4 helicopters only increase instance feasibility for the case of $CD^{\max}_i = 1$. Thus, from a feasibility point of view, the usage of two helicopters appears sufficient if expecting beneficiaries to walk up to 2 hours to pick up their aid is acceptable and the number of beneficaries does not exceed prior estimations. \\
While the introduction of a third or fourth helicopter only provides a marginal benefit for instance feasibility, they can allow the operator to reduce covering distances when $\rho_i$ is sufficiently small. In the following sections, we further elaborate on this circumstance.

\subsection{Impact of Additional Helicopters}\label{sec:additional_heli_impact}
Fleet sizing is a particularly crucial topic in emergency response, where budgets are often tight. While a small fleet incurs small operational costs, it leads to the need to aggregate demand toward fewer landing sites, incurring larger covering distances and thus, larger burdens put on the beneficiaries. On the contrary, larger fleet sizes are more expensive to charter and operate. Thus, in the problem at hand, the introduction of additional helicopters is only beneficial if covering durations are weighted accordingly. To be precise, for large route cost scaling factors, covering distances only play a secondary role and the main focus is to minimize operational costs. Thus, an increase of fleet size only impacts the solutions of instances with small route cost scaling factors. For this purpose, in this subsection, we only consider instances with $\rho_i = 0.1$.\\
\begin{figure}[ht]
\centering
\includegraphics[scale=0.55]{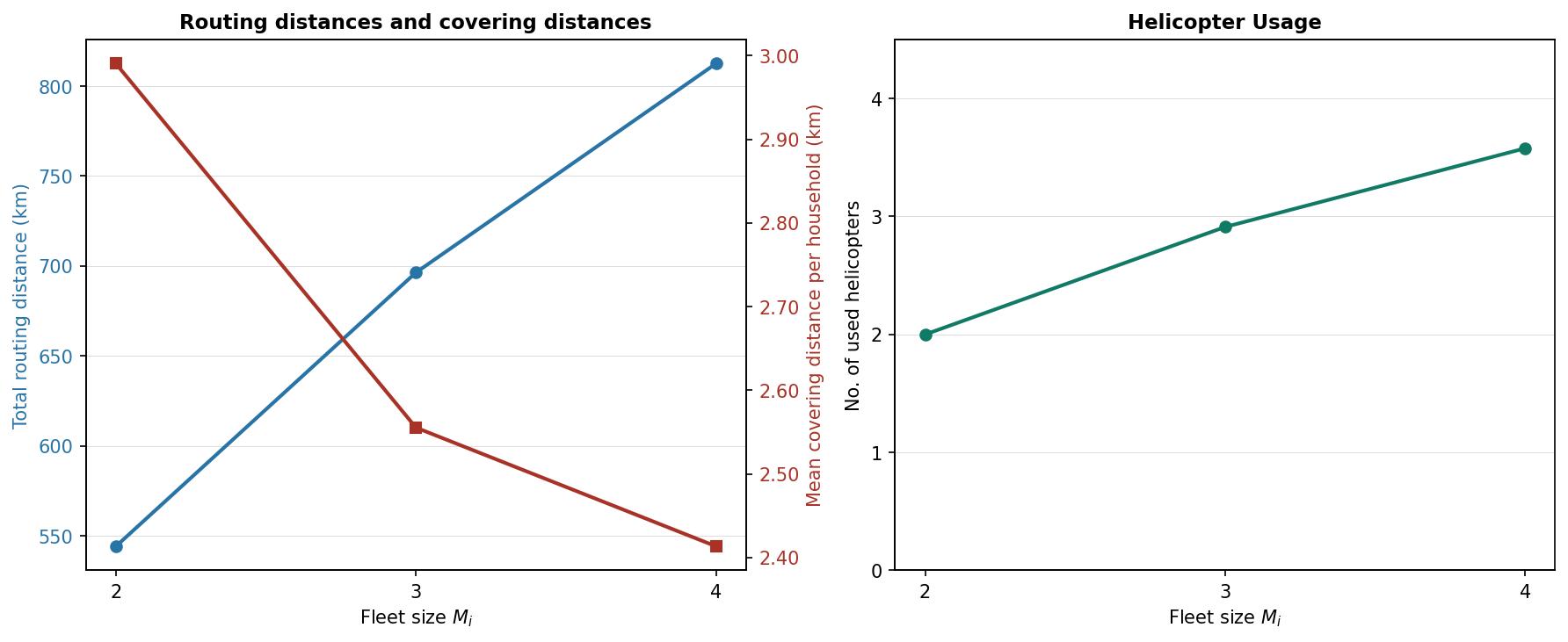} 
\caption{Distances and fleet utilization for different fleet sizes and $\rho_i = 0.1$}
\label{fig:vcount_analysis_1}
\end{figure}
\Cref{fig:vcount_analysis_1} visualizes the change in routing distances and mean covering distances per household, as well as the number of actually used helicopters for $\rho_i = 0.1$ and a fleet of 2, 3, and 4 available helicopters, respectively. 
To ensure comparability of results across different fleet sizes, \Cref{fig:vcount_analysis_1} only considers instances with $M_i \geq 2$ and $CD^{\max}_i \geq 2$, as for those parameter choices, for all instances with arbitrary $N^{\text{ben}}_i$ and $ST^{\max}_i$ a solution was found and thus the underlying instance sets have the same cardinality and only differ in the number of available helicopters. An increase of the fleet size from 2 to 3 helicopters reduces the mean covering distance per household by around 450 m, i.e., 15\%, and increases the routing distance by roughly 150 km, i.e., 27\%. Moreover, almost all solutions fully utilize the entire fleet (green line in the right-hand-side plot). On the contrary, when adding a fourth helicopter to the available fleet, covering distances only decrease by 150 m (i.e., by 6\%), while total routing distances increase by around 120 km (i.e., by 17\%). At the same time, the average number of used helicopters is around 3.6, indicating that a fourth helicopter is sometimes not beneficial, even when covering costs significantly dominate routing costs in the objective function. In conclusion, when covering efficiency plays an important role during devising the response strategy, the usage of 3 helicopters is an alternative worth considering. A fourth helicopter does not appear to bring a significant benefit, and its associated costs most likely dominate its benefit for covering distances, thus making it a less viable choice. Finally, we note that similar results can be observed for $\rho_i \in\{1, 3\}$. For $\rho_i > 3$, i.e., for $\rho\in\{5,7\},$ differences in routing distances, covering distances and the number of used helicopters become marginal, because additional helicopters are almost never utilized. Generally, the larger $\rho_i$ gets, the smaller the changes in routing and covering distances get, and the smaller the average number of used helicopters get.

\subsection{Trade-off between Operational Costs and Beneficiaries’ Burden}\label{sec:tradeoff_routing_vs_covering}
In order to find a fair balance between the operator's costs and the aforementioned burdens on the beneficiaries, in the following the impact of the previously introduced route cost scaling factor $\rho_i$ is evaluated.\\
\begin{figure}[ht]
\centering
\includegraphics[scale=0.55]{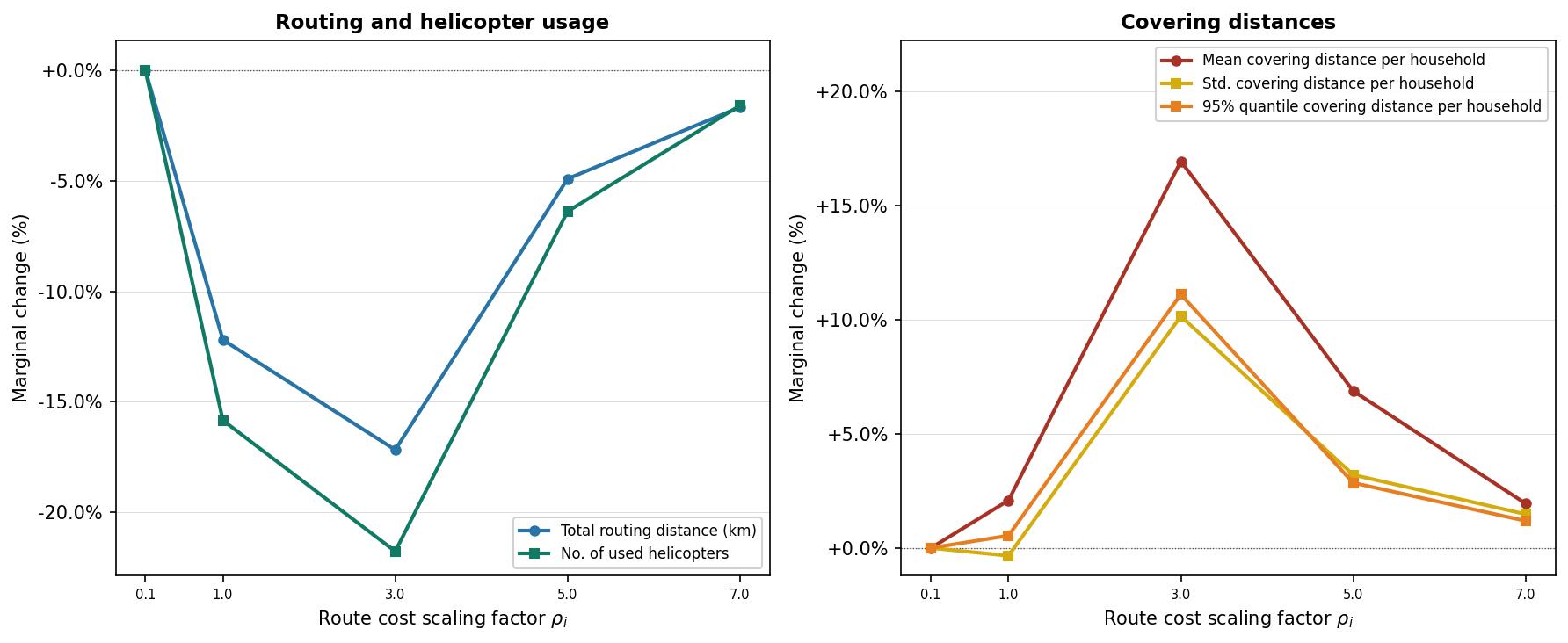} 
\caption{Marginal Changes of routing and covering values}
\label{fig:marginal_changes_routing_covering}
\end{figure}
\Cref{fig:marginal_changes_routing_covering} provides a summary of the marginal changes of several performance values when the route cost scaling factor $\rho_i$ is increased. Here, marginal changes refer to changes in the reference value when $\rho_i$ is increased by one step. For instance, the total routing distance decreases by an average of around 17\% when $\rho_i$ is increased from 1 to 3.\\
It can be seen that an increase in $\rho_i$ puts additional emphasis on few, short routes and reduces the relevance of covering decisions, leading to decreasing routing distances and number of used helicopters, while simultaneously increasing the mean and standard deviation of covering distances. In other words, when $\rho_i$ increases, routing distances are substituted with covering distances.\\
The right-hand side plot of \Cref{fig:marginal_changes_routing_covering} also allows for an analysis of the relationship between $\rho_i$ and several fairness-related values, such as the mean, the standard deviation and the  95\%-quantile of covering distances per household. It is clear to see that an increase of $\rho_i$, and thus an increasing focus on reducing operational costs, comes not just at the cost of increased covering distances for all beneficiaries, but it also increases dispersion between individuals. For instance, when increasing $\rho_i$ from 1 to 3, the average individual needs to walk around 17.5\% longer to reach their assigned landing site, where they can pick up their aid. At the same time, the standard deviation of covering distances as well as the 95\%-quantile of covering distances rise accordingly. This indicates that an increasing focus on routing costs does not impact all beneficiaries equally, but some beneficiaries are significantly more affected than others. Hence, when deciding on the balance between operational costs and burdens put on the beneficiaries, the operator must take into account not just the increase in average covering distances, but also the reduction of fairness among the beneficiaries.\\
\begin{figure}[ht]
\centering
\includegraphics[scale=0.45]{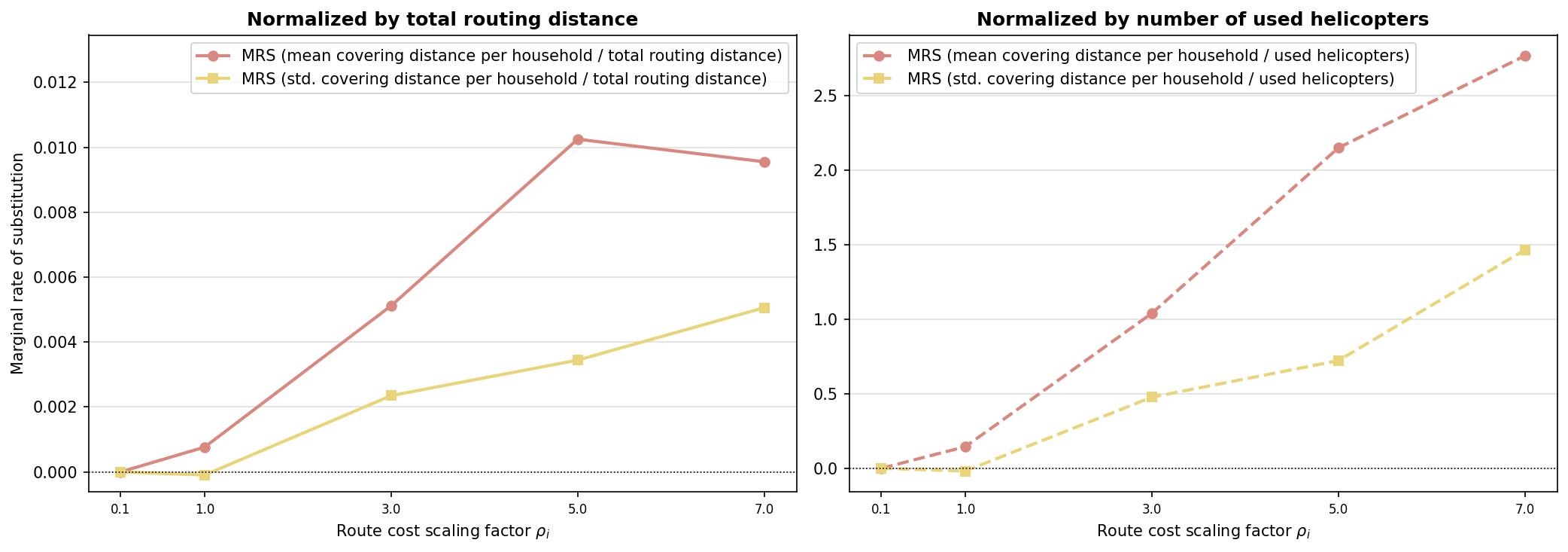} 
\caption{Marginal rates of substitution (MRS)}
\label{fig:rcs_mrs}
\end{figure}
 \Cref{fig:rcs_mrs} visualizes several discrete marginal rates of substitution (MRSs). At the example of the red curve on the left hand side plot, the MRS is computed as
 \begin{equation*}
     \text{MRS}(\rho_i) \coloneqq \dfrac{\text{mcd}(\rho_i) - \text{mcd}(\underline{\rho_i})}{\text{rd}(\underline{\rho_i}) - \text{rd}(\rho_i)}, 
 \end{equation*}
where $\text{mcd}(\rho_i)$ and $\text{rd}(\rho_i)$ are the mean covering distance per household and total routing distance for instances with route cost scaling factor $\rho_i$, and $\underline{\rho_i}$ is the largest evaluated route cost scaling factor smaller than $\rho_i$ (e.g., for $\rho_i = 7$, we obtain $\underline{\rho_i} = 5$). The MRS can be interpreted as the increase in covering distances one obtains per unit of routing distance (or used helicopter count) decrease when moving from $\underline{\rho_i}$ to $\rho_i$. In most cases, the MRS graphs are monotonically increasing, i.e., the amount of covering distance increase one needs to pay in order to reduce routing distances or the number of used helicopters gets larger for an increasing routing cost scaling factor. While the mean covering distance-based MRS is rather flat for $\rho\in\{0.1, 1\}$ and levels out for $\rho_i \geq 5$ when basing it on the total routing distance (left plot of \Cref{fig:rcs_mrs}), this is not the case for the standard deviation of covering distances, as well as for both the mean and standard deviation when the number of used helicopters is used as denominator (right plot of \Cref{fig:rcs_mrs}). 
This indicates that route cost scaling factors can be separated into three distinct regimes. In the first regime (i.e., $\rho \in [0,1]$), a decrease in operational costs is comparably cheap regarding the implied increasing ground travel burden. In the second regime (i.e., $\rho\in (1,5)$), operational costs can be further decreased, but covering distances start to increase more drastically. In the third regime (i.e., $\rho \geq 5$), reducing routing distances can be achieved by essentially paying a constant fee in the form of covering distances, whereas reducing the number of used vehicles becomes significantly more expensive. This indicates that, for practitioners who aim to find a suitable trade-off between operational costs and burden on beneficiaries, selecting a scaling factor $\rho\in [1,5]$ appears to be most appropriate. In either case, the proposed concept of scaling routing costs to achieve different trade-offs between operational costs and burdens on beneficiaries enables the user to evaluate different prioritization strategies and better understand the impact of their decisions on ground travel burdens and fairness.

\subsection{Impact of Service Threshold: Assessing the Value of Selective Coverage}\label{sec:triage}
A common approach in resource-constrained humanitarian operations is to define an operational coverage zone, beyond which beneficiaries are excluded from the response plan. In the model at hand, this is governed by the service threshold $ST^{\text{max}}_i$. Deciding whether such an approach is beneficial and to which extent it should be performed is often a difficult decision. Before considering it, its potential benefits and drawbacks must be analyzed.\\
\begin{figure}[ht]
\centering
\includegraphics[scale=0.6]{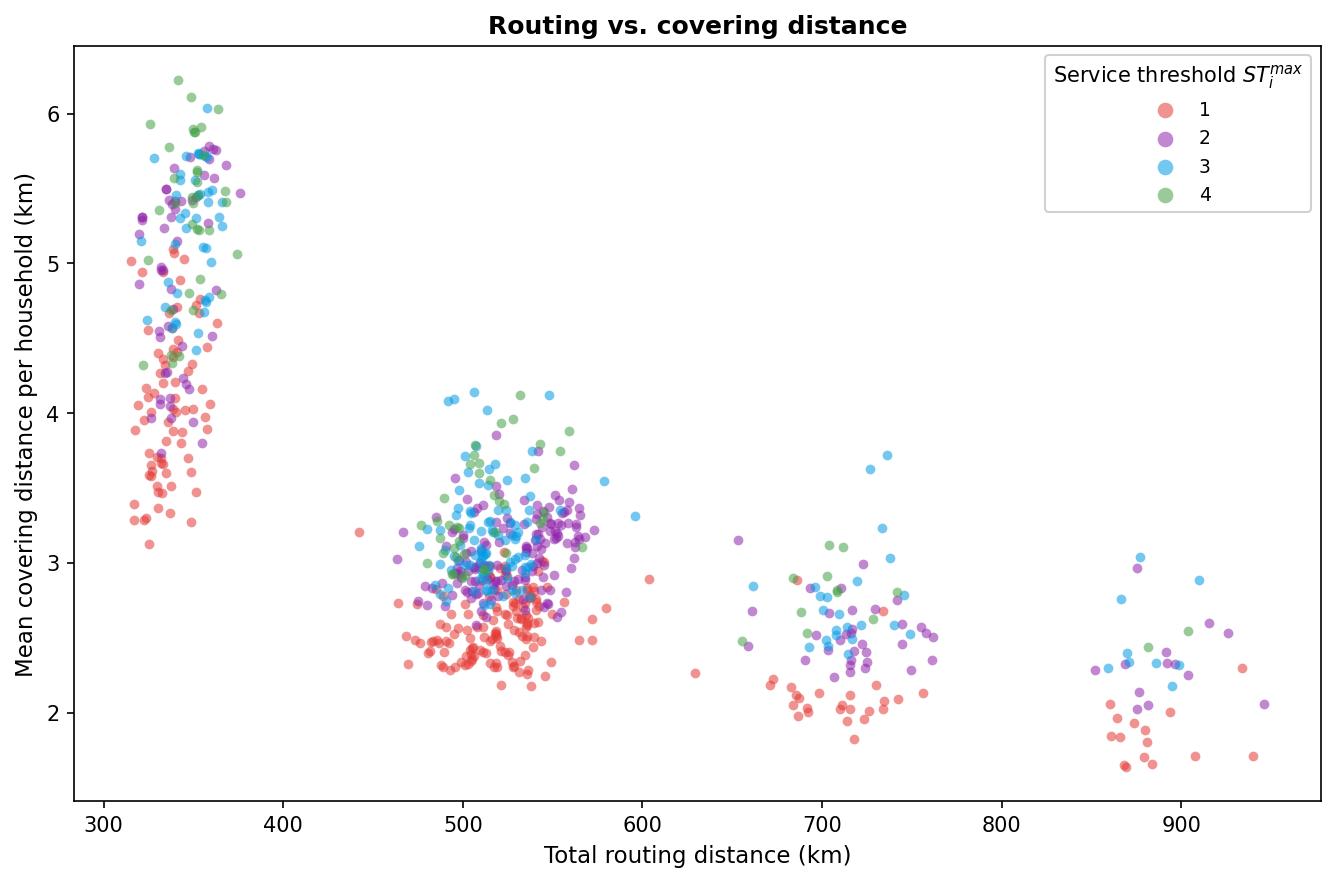} 
\caption{Scatter plot colored by service threshold}
\label{fig:scatter_plot_by_service_threshold}
\end{figure}
\Cref{fig:scatter_plot_by_service_threshold} shows a scatter plot of all instances with respect to total routing distances and mean covering distances per household of the best solution found. Each instance for which a solution was found is represented by one dot, and the dot colors reflect the underlying parameter value for the service threshold $ST^{\max}_i$. It can be seen that an increase of the service threshold shifts the scatter toward the top right corner. In other words, an increase of the service threshold increases the mean covering distance per household noticeably, while total routing distances increase slightly. This effect is particularly visible for $ST^{\max}_i\in \{1,2,3\}$. Note that instances with $ST^{\max}_i=3$ and $ST^{\max}_i=4$ only differ by 4 households. As stated at the beginning of this section, varying $ST^{\max}$ between 1 and 4 increases the number of beneficiaries from 5369 to 5831, i.e., around 8.6\%.\\
\begin{figure}[ht]
\centering
\includegraphics[scale=0.55]{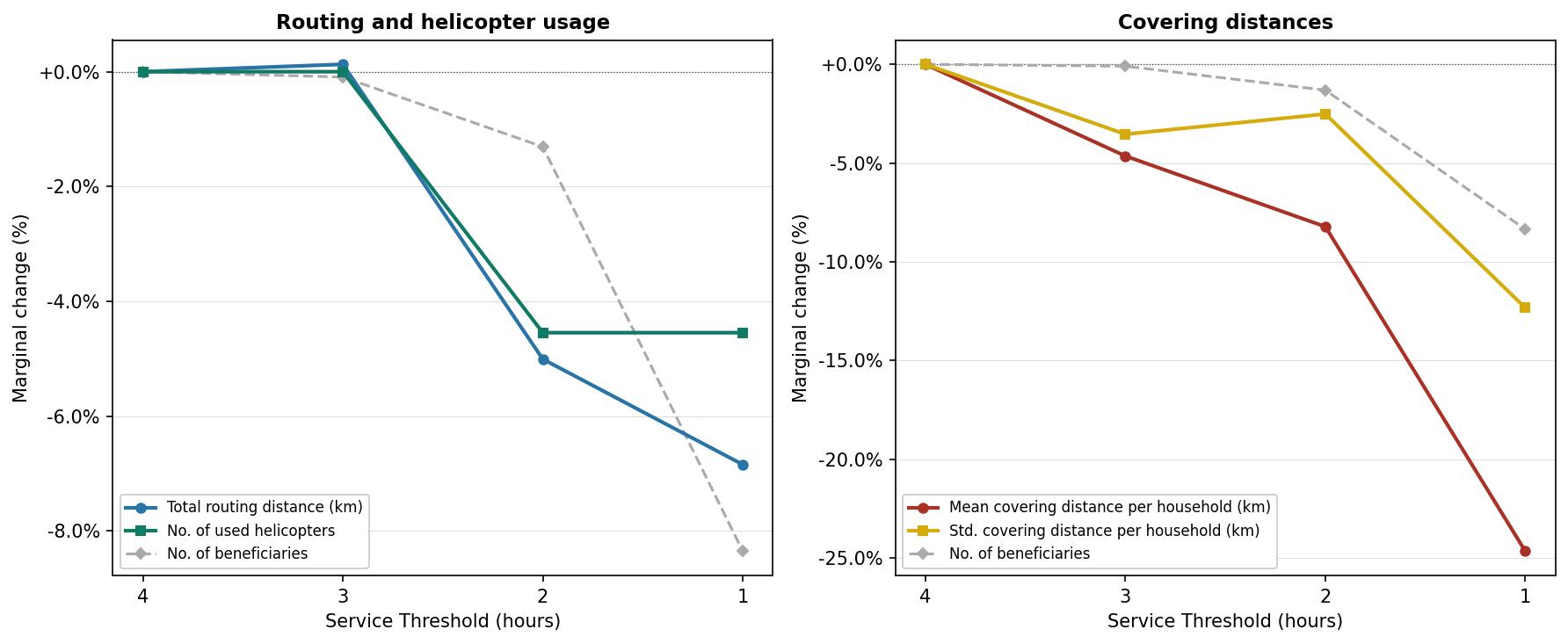} 
\caption{Marginal changes for $CD^{\max}_i = 4$ and $\rho_i\in\{5,7\}$}
\label{fig:margin_changes_cd4}
\end{figure}
\Cref{fig:margin_changes_cd4} provides a summary of the marginal changes of several performance values when the service threshold $ST^{\max}_i$ is decreased. Note that, because the introduction of a service threshold is typically only considered when reducing cost is of utmost importance, we only consider instances with $\rho_i \in \{5,7\}$, as \Cref{sec:tradeoff_routing_vs_covering} has shown that such route cost scaling factor choices put particular emphasis on routing cost. Furthermore, for the sake of simplicity, we only consider instances with $CD^{\max}_i = 4$. However, the following results also hold for instances with $CD^{\text{max}}_i < 4$. Similar to \Cref{fig:vcount_analysis_1}, \Cref{fig:margin_changes_cd4} only considers instances with $M_i \geq 2$ to ensure a stable baseline for comparison. Recall that all affected individuals have at least one landing site located within a 4 hour radius of their location. Taking the decrease in the number of modeled individuals (grey lines) into account, introducing a service threshold only improves efficiency (measured by routing distance [or no. of used helicopters] per beneficiary) for $ST^{\max}_i \geq 2$. At the same time, a decreasing service threshold decreases both the mean as well as the standard deviation of covering distances per household. Jointly with the results derived from \Cref{fig:scatter_plot_by_service_threshold}, one can conclude that imposing a service threshold of 2 hours reduces total operational costs and improves efficiency and ground travel burdens. However, this improvement comes at the price of disregarding the need for aid for a share of the affected population. Because the latter does not follow the core principles of humanitarian aid, introducing a service threshold should only be considered when minimizing cost is of utmost necessity.

%% file: sections/conclusion.tex
\section{Conclusion}\label{sec:conclusion}
In this paper, we studied a SCRP variant motivated by post-disaster humanitarian logistics, combining features of the m-CTP and the VRDAP in a multi-depot and facility-capacitated setting. To address the computational complexity of the problem, we developed an adaptive variable neighborhood search that integrates tabu search, simulated annealing and adaptive operator selection mechanisms with both established and novel local search operators, including a facility string replacement operator specifically designed for joint routing and covering decisions.
Benchmarking experiments demonstrated the competitiveness of the proposed AVNS across multiple problem classes. On m-CTP and m-CTP-p instances, the algorithm achieved near-optimal performance within a fraction of the runtime required by state-of-the-art exact approaches. On VRDAP instances, it identified 34 new best-known solutions. Notably, the algorithm found feasible solutions for 66 m-CTP-p instances for which no feasible solution was previously known. These results demonstrate the robustness and broad applicability of the proposed approach across structurally different problem variants.\\
Furthermore, we applied the proposed algorithm to a case study in the 2024 flash flood events in central Afghanistan. We utilize a non-monetary weighted objective function that enables the modeling of different trade-offs between operational costs and beneficiary burden. The analysis provides several relevant insights for practitioners, including the observation that a fleet of two helicopters is generally sufficient for feasible coverage when beneficiaries can walk for up to 2 hours, while a third helicopter offers meaningful reductions in covering distances. We also conclude that a service threshold of 2 hours can improve operational efficiency but at the cost of excluding affected individuals, and should thus be used only when necessary. Overall, the results demonstrate the potential of the proposed framework as a flexible decision-support approach for humanitarian logistics settings characterized by limited accessibility and multiple interacting decision layers.\\
From an application point of view, future research directions include incorporating uncertain demand and beneficiary locations that evolve over time. Furthermore, extending the model to multi-period planning horizons and investigating prepositioning strategies, in combination with the proposed helicopter aid approach, to further reduce response times in anticipation of recurring flood events can be of particular interest. From a methodological perspective, an extension of the proposed framework to additional SCRP variants and a comparison to state-of-the-art approaches for these problems can be of particular interest.

%% file: sections/acknowledgement.tex
\section*{Acknowledgements}
\noindent Andreas Hagn has been funded by the Deutsche Forschungsgemeinschaft (DFG, German Research Foundation) - Project number 277991500. Jan Krause and Moritz Stargalla were funded by the Bayerische Forschungsstiftung within the SHIELD project. \\

%% file: sections/supplementary.tex
\section{List of Symbols} \label{sec:notation}
\begin{longtable}{ll}

\endfirsthead

\endhead

\endfoot

\endlastfoot

\small  
\text{ }\\
\multicolumn{2}{@{}l}{\textbf{Sets}} \\
$D$ & Set of depots\\
$F$ & Set of facilities\\
$W$ & Set of customers\\
$F_l$ & Set of facilities that can cover customer $l\in W$\\
$A$ & Set of arcs connecting depots and facilities\\
$(\mathcal{R}, \phi)$ & Feasible solution, given by a set $\mathcal{R}$ of routes and a customer assignment $\phi$\\
$\mathcal{D}$ & Set of selected depots (implied by a solution $(\mathcal{R}, \phi)$\\
$\mathcal{F}$ & Set of selected facilities (implied by a solution $(\mathcal{R}, \phi)$\\
\\
\multicolumn{2}{l}{\textbf{Variables}} \\
$x_{ij}$ & Binary arc selection variables\\
$z_{li}$ & Binary customer assignment variables\\
$w_i$ & Binary facility selection variables\\
$u_i$ & Load tracking variables\\
$\pi_i$ & Route stops tracking variables\\
$\delta_i$ & Route length variables \\
\\
\multicolumn{2}{l}{\textbf{Model Parameters}} \\
$c_{ij}$ & Edge traversal cost\\
$\text{dist}_{i,j}$ & Distance between nodes $i$ and $j$\\
$r_l$ & Customer demand\\
$Q$ & Vehicle capacity\\
$Q_i^F$ & Facility capacity of facility $i\in F$\\
$M$ & No. of available vehicles\\
$p$ & Maximum number of stops (i.e., facilities) per route\\
$q$ & Maximum length (distance-wise for case study instances, routing cost-wise else) per route\\
\\
\multicolumn{2}{l}{\textbf{AVNS Parameters}} \\
$N^{\text{start}}$ & Maximum number of constructive heuristic attempts\\
$\mathcal{O}$ & Set of available operators\\
$\lambda^{\text{initial}}$ & Initial score assigned to each operator\\
$\lambda = (\lambda_o)_{o\in\mathcal{O}}$ & Vector of current scores of all operators\\
$\sigma$ & Score update vector\\
$\lambda^{\text{seg}}$ & Vector of operator scores in the current segment\\
$N^{\text{seg}}$ & Segment length (no. of iterations after which operator scores $\lambda$ are updated)\\
$\rho$ & Smoothing rate (for score updates at segment end)\\
$n_{\max}$ & Maximum number of accepted changes per facility per FacilityStringReplacement call\\
$N^{\text{stag}}$ & Maximum number of iterations without solution improvement until shaking is called\\
$N^{\text{shaking}}$ & Maximum number of shaking attempts per shaking call\\
$(N_k)_{k=1,\dots,K}$ & Shaking neighborhoods\\
$T^{\text{initial}}$ & Initial temperature\\
$T^{\text{start}}$ & Start temperature (obtained by multiplying $T^{\text{initial}}$ with initial solution cost)\\
$T$ & Current temperature\\
$\mathcal{L}$ & Tabu list\\
$N^{\text{tabu}}$ & Tabu list size\\
$t^{\text{elap}}$ & Elapsed runtime (in CPU seconds)\\
$t^{\text{max}}$ & Maximum runtime (in CPU seconds)\\
$\left(\mathcal{R}^{\text{best}}, \phi^{\text{best}}\right)$ & Best solution found by the AVNS\\

\end{longtable}

\clearpage

\section{Case Study: Instance Generation}
\subsection{Instance Generation}

In the following, we provide additional details on the logic that was used to construct test instances for the case study presented in Section 6 of the manuscript.\\
\paragraph{\textbf{Fleet Size and Route Lengths}}\text{ }\\
Typically, the World Food Programme (WFP) uses aircraft operated by the UN Humanitarian Air Service (UNHAS) to carry out airborne logistics. As a result, until 2024, the UNHAS operated an Mil Mi-8 helicopter in Afghanistan (cf. \cite{WFP2024UNHAS}). Furthermore, variants of said helicopter type are particularly suited for missions in rough circumstances, which makes them a good fit for operations in the mountainous, high-elevation areas of the considered districts in Baghlan province. Therefore, we assume a homogeneous fleet of 1 to 4 Mi-8AMT helicopters with a range of 540 km (cf. \cite{heli_datasheet}). Hence, we set the maximum route length to be $q\coloneqq 540$. According to \cite{heli_datasheet}, the maximum payload of a Mil Mi-8 AMT is 4,000 kg. However, helicopter performance tends to degrade in higher elevation areas due to reduced air density. To ensure operability, we thus reduce the payload by 20\%, yielding a payload of 3,200 kg.\\
The maximum number of stops $p$ is derived as follows. We consider a maximum operating time window of 8 hours per helicopter per day. Assuming an average travel speed of 170 $km/h$ and a range of $q=570$ km, the maximum flying time, excluding time for arrival at and departure from landing sites, per day is $\dfrac{570}{170}=3.35$ hours, or 201 minutes. This leaves around 4.65 hours (or 279 minutes) for landing, unloading, and departing. We assume that landing at a landing site, unloading aid, and departing takes 10 minutes on average. This includes time required for descent, landing, unloading, departure, and ascent. Thus, the maximum number of stops a helicopter can make considering the operating time of up to 8 hours is $\dfrac{279}{10} = 27.9 \approx 28$.\\

\paragraph{\textbf{Customer Locations and Demands}}\text{ }\\
We obtained data on the number of affected individuals per district during the May 2024 Floods, as well as population data for areas of size 1 km² in Afghanistan from \cite{affected_individuals} and \cite{population_data}, respectively. For each district, we calculate the ratio of affected individuals versus its population. Then, for each area of size 1 km², the number of affected individuals (i.e. beneficiaries) is given by multiplying said ratio with its population. Each such area with at least one affected individual then corresponds to a customer node. Depending on the instance, the number of beneficiaries is then multiplied by $N^{\text{ben}}_i > 0$.\\
The demand of a customer node is calculated as the demand per individual multiplied by the number of affected individuals. The demand per person is calculated as the weight of high-energy biscuits (HEBs) required to cover their minimal daily basic calorie intake. The \cite{min_calorie_intake} specify said intake to be 2,100 kcal, while high-energy biscuits have a energy density of 450 kcal per 100 g (cf. \cite{heb_density}). We note that we restrict the aid distribution to HEBs because in the considered scenario, the UN's immediate response also relied heavily on HEBs (cf. \cite{wfp_heb_response}). Therefore, each individual is assumed to have a demand of $\dfrac{2100}{450}\cdot 0.1 \approx 0.467$ kg. Overall, instances have between 432 and 558 customer nodes, representing 2819 to 8922 beneficiaries with a total daily demand of 1316.47 kg to 4166.57 kg.\\

\paragraph{\textbf{Facility Locations}}\text{ }\\
We scan the entire affected region of Baghlan and its surroundings for suitable landing sites. We assume that a landing site needs to have a minimum size of approximately 100 m x 100 m, i.e. 3 arcseconds x 3 arcseconds. Furthermore, landing sites must be located on grasslands or shrublands, with no major vegetation that might hinder safe landing or take-off. Each landing site must have a minimum distance of 150 m to waterbodies, as these are often overflown during flash floods, as well as a minimum distance of 70 m to any other terrain that makes helicopter approaches dangerous, such as mountains or areas with dense vegetation. Furthermore, the local curvature around the potential landing site must not be below -0.005, as such areas are typically valleys prone to flooding. We note that -0.005 is the 90\%-quantile of the distribution of local curvatures in the affected region and thus a reasonable proxy for the identification of flood-prone valleys. Finally, we also require landing sites to have a maximum local slope of 7\% and to not be located above 5,400 m, as Mil Mi-8 AMT helicopters can only safely operate up to said elevation assuming maximum payload (cf. \cite{heli_datasheet}). To evaluate the aforementioned criteria, \cite{tandem_elevation_data} and \cite{esa_terrain_data} provide publicly available data on elevations and terrain types, respectively, across Afghanistan. We note that said elevation data is collected via satellite imagery, which, at times, is rather prone to imprecisions. This frequently leads to significant `jumps' in elevation in small areas, even though terrain data does not indicate mountainous terrain. An example for this circumstance are meadows with large trees, where elevation computations are often disturbed by treetops, leading to large elevation changes in small areas. We take this into account by dropping the aforementioned maximum slope criterion if the slope within a landing site exceeds 30\%. Note that the usage of a terrain filter as described before mitigates the risk of identifying mountainous regions as suitable landing sites.
\\
Moreover, we assume that last-mile distribution from the landing sites to the affected individuals still relies on road networks. Even though the affected regions are inaccessible from the outside, transport within them can be carried out on roads by foot or on mules at a reduced speed of 4 km/h. Since roads may be partially impassable and circumvention is required, we increase the incurred ground travel distances by 20\% in the following. Furthermore, because landing sites typically do not have a paved access to the nearest road, we impose that the haversine distance between any landing site and the nearest road must not exceed 1,000 m and does not cross any inaccessible terrain, such as mountains or waterbodies. For this purpose, we collected data on the Afghan road network from \cite{road_networks}. Because the resulting graph is highly disconnected and contains many isolated edges, it is post-processed by iteratively identifying pairs of connected components and linking them via the shortest (artificial) edge that neither crosses forbidden terrain (waterbodies, mountains, areas with dense vegetation) nor exceeds a maximum slope of 15 degrees.\\
Each landing site can then be mapped to its closest node on the road network. In the following, we call the said road network nodes to which at least one potential landing site is mapped \textit{facility candidates}. Note that this implies the existance of a surjective mapping from landing sites to facility candidates, i.e., each facility candidate can equivalently be represented by at least one landing site. After applying all of the aforementioned filters, we obtain 20,430 facility candidates. Because many of these candidates are very close to each other, a large majority of them can be safely removed without significantly compromising solution quality. Thus, we select a subset of up to 500 facility candidates using the KMeans algorithm in the \textit{scikit-learn} Python package. The feature vector consists of the distances between the facility candidate and all customer nodes. If any customer nodes are not covered within the maximum covering duration by any of the selected facility candidates, we iteratively add the closest candidate for each not-yet-covered customer until every customer is covered by at least one facility. Because many landing sites are often in close proximity to each other and do not fundamentally differ with respect to the beneficiaries they can cover, we then remove all landing sites (or, more precisely, their corresponding facility candidate) for which another landing site (that corresponds to a different facility candidate) is within a 12-minute walking radius. Distribution on-site is assumed to be undertaken by the beneficiaries themselves, where we assume that out of each household, one person picks up their aid. In accordance with \cite{afg_household_data}, a household in Baghlan province is assumed to consist of 6 individuals. To limit the (unloading) work required on-site and to avoid single points of failure in the distributional structure, we limit the capacity of each landing site to $Q^F_f = \left\{\underset{l\in W: \ f\in F_l}{\sum}r_l,200\right\}$ kg for each facility candidate $f$, where the first sum equals the total demand of all customer nodes coverable by the facility candidate.

\paragraph{\textbf{Depot Locations}}\text{ }\\
Depot locations correspond to international, major, regional, or local airports or helipads that have a haversine distance of at most 150 km to at least one affected individual. A total of 14 locations satisfy this condition. We assume that depots are uncapacitated and their usage comes at no cost. Thus, an optimal solution never contains a depot that is not the closest depot to any facility candidate. Removing these guaranteed-to-be-suboptimal locations yields 2 potential depots: Taloqan Airport and Bagram Airfield.

\paragraph{\textbf{Routing and Covering Costs}}\text{ }\\
In the considered situation, the operator needs to balance the operational costs (consisting of distance-based costs, such as fuel, and fixed costs, such as crew salary and helicopter chartering) and the burdens put on the beneficiaries by requiring them to travel to their designated landing site in order to pick up their aid. Hence, we assume that both routing aircraft between landing sites and assigning beneficiaries to landing sites incur cost. However, there are 2 things to consider. First, information on monetary cost of helicopter chartering, fuel, or airfield usage in Afghanistan is hardly available and unreliable. Second, while operating costs can be measured in monetary terms, burdens caused to beneficiaries cannot. Thus, the proposed instances use routing and covering costs that are designed to be proxies for the aforementioned cost dimensions. In order to so, these costs need to satisfy several properties: (a) the cost of an additional helicopter always exceeds the potential savings realized by splitting an existing route into two routes (i.e., fewer routes are always cheaper), (b) helicopter travelling costs are directly proportional to the travelled distance, and (c) covering costs are proportional to the caused burdens, measured by ground travel times and the number of beneficiaries that need to undertake ground travel. Therefore, we consider the following routing costs:
\begin{align*}
    c_{ij} = \begin{cases}
        \begin{aligned}
            &\rho_i \cdot \text{dist}_{i,j} && \text{if} \ \ i,j\notin D\\
            &\rho_i \cdot \left(\text{dist}_{i,j} + \dfrac{1}{2}\cdot\left(\text{dist}^F_{\max} + 2\cdot\text{dist}^D_{\max}\right)\right) && \text{else}
        \end{aligned}
    \end{cases}
\end{align*}
where $\rho_i \geq 0 $ is an instance-dependent rescaling factor called \textit{route cost scaling factor}, $\text{dist}_{i,j} \geq 0$ is the helicopter travel distance between facilities (or depots) $i,j\in D\cup F$, and
\begin{align*}
    \text{dist}^F_{\max} &\coloneqq \max\left\{\text{dist}_{i,j}: i,j\in F\right\}\\
    \text{dist}^D_{\max} &\coloneqq \max\left\{\text{dist}_{d,i}: d\in D, i\in F\right\}
\end{align*}
are the maximum facility-to-facility- and depot-to-facility-distances. The latter constant term for the case that $i\in D$ or $j\in D$ holds acts as a fixed cost for selecting a vehicle. It is designed to ensure that splitting an existing route into two routes always increases routing costs. To understand its functionality, consider the case of splitting a route into two new routes and potentially choosing new depots to start the routes from. By doing so, one facility-to-facility edge and two depot-to-facility edges are removed from the solution, while four new depot-to-facility edges are introduced. Therefore, it is easy to see that the potential cost reduction incurred by this route split is bounded from above by $\text{dist}^F_{\max} + 2 \cdot \text{dist}^D_{\max}$ and thus does not exceed the cost of adding a new route.\\
Furthermore, covering costs are computed as
\begin{align*}
    a_{li} \coloneqq \text{dist}_{l,i} \cdot \left\lceil \dfrac{B_l}{H}\right\rceil
\end{align*}
where $B_l \in\mathbb{N}$ is the number of beneficiaries represented by customer node $l\in W$ and $H = 6$ is the household size.\\
Note that $\rho_i$ is instance-dependant and ranges between $0.1$ and $7$. Preliminary studies have shown that instances with $\rho_i = 0.1$ (respectively $\rho_i = 7$) yield solutions that exhibit total routing and covering distances similar to the ones of instances, where routing (respectively covering) costs are set to zero. Hence, these values for $\rho_i$ represent two extreme valuation schemes, where one of the two objective dimensions is almost entirely ignored.\\

%% file: bibliography.bib
@article{Oliveira2025_mCTP_exact,
  title = {New cuts and a branch-cut-and-price model for the multi-vehicle covering tour problem},
  ISSN = {1614-2411},
  DOI = {10.1007/s10288-025-00584-0},
  journal = {4OR},
  publisher = {Springer Science and Business Media LLC},
  author = {Oliveira,  Bruno and Pessoa,  Artur and Roboredo,  Marcos},
  year = {2025},
  month = jan 
}

@article{Pham2017,
  title = {Solving the multi-vehicle multi-covering tour problem},
  volume = {88},
  ISSN = {0305-0548},
  DOI = {10.1016/j.cor.2017.07.009},
  journal = {Computers \& Operations Research},
  publisher = {Elsevier BV},
  author = {Pham,  Tuan Anh and Hoàng Hà,  Minh and Nguyen,  Xuan Hoai},
  year = {2017},
  month = dec,
  pages = {258–278}
}

@article{Ghoniem2013,
  title = {A specialized column generation approach for a vehicle routing problem with demand allocation},
  volume = {64},
  ISSN = {1476-9360},
  DOI = {10.1057/jors.2012.32},
  number = {1},
  journal = {Journal of the Operational Research Society},
  publisher = {Informa UK Limited},
  author = {Ghoniem,  A and Scherrer,  C R and Solak,  S},
  year = {2013},
  month = jan,
  pages = {114–124}
}

@article{AnayaArenas2014,
  title = {Relief distribution networks: a systematic review},
  volume = {223},
  ISSN = {1572-9338},
  DOI = {10.1007/s10479-014-1581-y},
  number = {1},
  journal = {Annals of Operations Research},
  publisher = {Springer Science and Business Media LLC},
  author = {Anaya-Arenas,  A. M. and Renaud,  J. and Ruiz,  A.},
  year = {2014},
  month = apr,
  pages = {53–79}
}

@article{Hoyos2015,
  title = {OR models with stochastic components in disaster operations management: A literature survey},
  volume = {82},
  ISSN = {0360-8352},
  DOI = {10.1016/j.cie.2014.11.025},
  journal = {Computers \& Industrial Engineering},
  publisher = {Elsevier BV},
  author = {Hoyos,  Maria Camila and Morales,  Ridley S. and Akhavan-Tabatabaei,  Raha},
  year = {2015},
  month = apr,
  pages = {183–197}
}

@article{Davoodi2019,
  title = {An integrated disaster relief model based on covering tour using hybrid Benders decomposition and variable neighborhood search: Application in the Iranian context},
  volume = {130},
  ISSN = {0360-8352},
  DOI = {10.1016/j.cie.2019.02.040},
  journal = {Computers \& Industrial Engineering},
  publisher = {Elsevier BV},
  author = {Davoodi,  Sayyed Mohammad Reza and Goli,  Alireza},
  year = {2019},
  month = apr,
  pages = {370–380}
}

@article{Sacramento2019,
  title = {An adaptive large neighborhood search metaheuristic for the vehicle routing problem with drones},
  volume = {102},
  ISSN = {0968-090X},
  DOI = {10.1016/j.trc.2019.02.018},
  journal = {Transportation Research Part C: Emerging Technologies},
  publisher = {Elsevier BV},
  author = {Sacramento,  David and Pisinger,  David and Ropke,  Stefan},
  year = {2019},
  month = may,
  pages = {289–315}
}

@article{Vidal2013,
  title = {Heuristics for multi-attribute vehicle routing problems: A survey and synthesis},
  volume = {231},
  ISSN = {0377-2217},
  DOI = {10.1016/j.ejor.2013.02.053},
  number = {1},
  journal = {European Journal of Operational Research},
  publisher = {Elsevier BV},
  author = {Vidal,  Thibaut and Crainic,  Teodor Gabriel and Gendreau,  Michel and Prins,  Christian},
  year = {2013},
  month = nov,
  pages = {1–21}
}

@incollection{Rasku2019,
  author    = {Joonas Rasku and Tommi K{\"a}rkk{\"a}inen and Nysret Musliu},
  title     = {Meta-Survey and Implementations of Classical Capacitated Vehicle Routing Heuristics with Reproduced Results},
  booktitle = {Toward Automatic Customization of Vehicle Routing Systems},
  pages     = {133--260},
  series    = {JYU Dissertations},
  volume    = {113},
  publisher = {University of Jyväskylä},
  address   = {Jyväskylä, Finland},
  year      = {2019}
}

@article{WindrasMara2022,
  title = {A survey of adaptive large neighborhood search algorithms and applications},
  volume = {146},
  ISSN = {0305-0548},
  DOI = {10.1016/j.cor.2022.105903},
  journal = {Computers \& Operations Research},
  publisher = {Elsevier BV},
  author = {Windras Mara,  Setyo Tri and Norcahyo,  Rachmadi and Jodiawan,  Panca and Lusiantoro,  Luluk and Rifai,  Achmad Pratama},
  year = {2022},
  month = oct,
  pages = {105903}
}

@book{bookmetaheuristics2019,
  editor = {Gendreau, M. and Potvin, J.-Y.},
  title = {Handbook of Metaheuristics},
  ISBN = {9783319910864},
  ISSN = {2214-7934},
  DOI = {10.1007/978-3-319-91086-4},
  journal = {International Series in Operations Research \& Management Science},
  publisher = {Springer International Publishing},
  year = {2019}
}

@article{Mladenovi1997,
  title = {Variable neighborhood search},
  volume = {24},
  ISSN = {0305-0548},
  DOI = {10.1016/s0305-0548(97)00031-2},
  number = {11},
  journal = {Computers \& Operations Research},
  publisher = {Elsevier BV},
  author = {Mladenović,  N. and Hansen,  P.},
  year = {1997},
  month = nov,
  pages = {1097–1100}
}

@misc{heuristicssurvey,
  doi = {10.48550/ARXIV.2303.04147},
  author = {Liu,  Fei and Lu,  Chengyu and Gui,  Lin and Zhang,  Qingfu and Tong,  Xialiang and Yuan,  Mingxuan},
  title = {Heuristics for Vehicle Routing Problem: A Survey and Recent Advances},
  publisher = {arXiv},
  year = {2023},
  copyright = {arXiv.org perpetual,  non-exclusive license}
}

@inbook{Nolz2010,
  title = {A Bi-objective Metaheuristic for Disaster Relief Operation Planning},
  ISBN = {9783642112188},
  ISSN = {1860-9503},
  DOI = {10.1007/978-3-642-11218-8_8},
  booktitle = {Advances in Multi-Objective Nature Inspired Computing},
  publisher = {Springer Berlin Heidelberg},
  author = {Nolz,  Pamela C. and Doerner,  Karl F. and Gutjahr,  Walter J. and Hartl,  Richard F.},
  year = {2010},
  pages = {167–187}
}

@article{Tricoire2012,
  title = {The bi-objective stochastic covering tour problem},
  volume = {39},
  ISSN = {0305-0548},
  DOI = {10.1016/j.cor.2011.09.009},
  number = {7},
  journal = {Computers \& Operations Research},
  publisher = {Elsevier BV},
  author = {Tricoire,  Fabien and Graf,  Alexandra and Gutjahr,  Walter J.},
  year = {2012},
  month = jul,
  pages = {1582–1592}
}

@article{Goli2019,
  title = {A Covering Tour Approach for Disaster Relief Locating and Routing with Fuzzy Demand},
  volume = {18},
  ISSN = {1868-8659},
  DOI = {10.1007/s13177-019-00185-2},
  number = {1},
  journal = {International Journal of Intelligent Transportation Systems Research},
  publisher = {Springer Science and Business Media LLC},
  author = {Goli,  Alireza and Malmir,  Behnam},
  year = {2019},
  month = may,
  pages = {140–152}
}

@article{Kl2025,
  title = {Modeling a humanitarian‐aid covering tour problem with location selection and vehicle assignment decisions},
  ISSN = {1475-3995},
  DOI = {10.1111/itor.70088},
  journal = {International Transactions in Operational Research},
  publisher = {Wiley},
  author = {Kılı\c{c},  Kaan and Meterelliyoz,  Melike and G\"{u}ven\c{c} Pelit,  Ilay and Soysal,  Mehmet},
  year = {2025},
  month = aug 
}

@article{Alinaghian2019,
  title = {A mathematical model for location of temporary relief centers and dynamic routing of aerial rescue vehicles},
  volume = {131},
  ISSN = {0360-8352},
  DOI = {10.1016/j.cie.2019.03.002},
  journal = {Computers \& Industrial Engineering},
  publisher = {Elsevier BV},
  author = {Alinaghian,  Mahdi and Aghaie,  Mohammad and Sabbagh,  Mohammad S.},
  year = {2019},
  month = may,
  pages = {227–241}
}

@article{delaTorre2012,
  title = {Disaster relief routing: Integrating research and practice},
  volume = {46},
  ISSN = {0038-0121},
  DOI = {10.1016/j.seps.2011.06.001},
  number = {1},
  journal = {Socio-Economic Planning Sciences},
  publisher = {Elsevier BV},
  author = {de la Torre,  Luis E. and Dolinskaya,  Irina S. and Smilowitz,  Karen R.},
  year = {2012},
  month = mar,
  pages = {88–97}
}

@article{Anuar2021,
  title = {Vehicle Routing Optimisation in Humanitarian Operations: A Survey on Modelling and Optimisation Approaches},
  volume = {11},
  ISSN = {2076-3417},
  DOI = {10.3390/app11020667},
  number = {2},
  journal = {Applied Sciences},
  publisher = {MDPI AG},
  author = {Anuar,  Wadi Khalid and Lee,  Lai Soon and Pickl,  Stefan and Seow,  Hsin-Vonn},
  year = {2021},
  month = jan,
  pages = {667}
}

@inproceedings{Xu2018,
  series = {cimns-18},
  title = {Emergency Logistics Theory,  Model and Method: A Review and Further Research Directions},
  DOI = {10.2991/cimns-18.2018.42},
  booktitle = {Proceedings of the 2018 3rd International Conference on Communications,  Information Management and Network Security (CIMNS 2018)},
  publisher = {Atlantis Press},
  author = {Xu,  Hongqian and Fang,  Danhui and Jin,  Yining},
  year = {2018},
  collection = {cimns-18}
}

@article{Zheng2015,
  title = {Evolutionary optimization for disaster relief operations: A survey},
  volume = {27},
  ISSN = {1568-4946},
  DOI = {10.1016/j.asoc.2014.09.041},
  journal = {Applied Soft Computing},
  publisher = {Elsevier BV},
  author = {Zheng,  Yu-Jun and Chen,  Sheng-Yong and Ling,  Hai-Feng},
  year = {2015},
  month = feb,
  pages = {553–566}
}

@article{Luo2023,
  title = {Optimization models and solving approaches in relief distribution concerning victims' satisfaction: A review},
  volume = {143},
  ISSN = {1568-4946},
  DOI = {10.1016/j.asoc.2023.110398},
  journal = {Applied Soft Computing},
  publisher = {Elsevier BV},
  author = {Luo,  Jia and Shi,  Lei and Xue,  Rui and El-baz,  Didier},
  year = {2023},
  month = aug,
  pages = {110398}
}

@article{Reihaneh2018,
  title = {A multi-start optimization-based heuristic for a food bank distribution problem},
  volume = {69},
  ISSN = {1476-9360},
  DOI = {10.1057/s41274-017-0220-9},
  number = {5},
  journal = {Journal of the Operational Research Society},
  publisher = {Informa UK Limited},
  author = {Reihaneh,  Mohammad and Ghoniem,  Ahmed},
  year = {2018},
  month = dec,
  pages = {691–706}
}

@article{Reihaneh2019,
  title = {A branch-and-price algorithm for a vehicle routing with demand allocation problem},
  volume = {272},
  ISSN = {0377-2217},
  DOI = {10.1016/j.ejor.2018.06.049},
  number = {2},
  journal = {European Journal of Operational Research},
  publisher = {Elsevier BV},
  author = {Reihaneh,  Mohammad and Ghoniem,  Ahmed},
  year = {2019},
  month = jan,
  pages = {523–538}
}

@article{Oliveira2025_VRDAP_exact,
  title = {An exact approach for the Vehicle Routing Problem with Demand Allocation},
  volume = {182},
  ISSN = {0305-0548},
  DOI = {10.1016/j.cor.2025.107101},
  journal = {Computers \& Operations Research},
  publisher = {Elsevier BV},
  author = {Oliveira,  Bruno and Lima,  Diogo and Pessoa,  Artur and Roboredo,  Marcos},
  year = {2025},
  month = oct,
  pages = {107101}
}

@article{Oliveira2025_VRDAP_ILS,
  title = {Hybrid iterated local search algorithm for the vehicle routing problem with lockers},
  volume = {31},
  ISSN = {1572-9397},
  DOI = {10.1007/s10732-025-09557-2},
  number = {2},
  journal = {Journal of Heuristics},
  publisher = {Springer Science and Business Media LLC},
  author = {Oliveira,  Bruno and Pessoa,  Artur and Roboredo,  Marcos},
  year = {2025},
  month = apr 
}

@article{Oliveira2015,
  title = {MULTI-VEHICLE COVERING TOUR PROBLEM: BUILDING ROUTES FOR URBAN PATROLLING},
  volume = {35},
  ISSN = {0101-7438},
  DOI = {10.1590/0101-7438.2015.035.03.0617},
  number = {3},
  journal = {Pesquisa Operacional},
  publisher = {FapUNIFESP (SciELO)},
  author = {Oliveira,  Washington Alves de and Moretti,  Antonio Carlos and Reis,  Ednei Felix},
  year = {2015},
  month = dec,
  pages = {617–644}
}

@article{MoshrefJavadi2016,
  title = {The Latency Location-Routing Problem},
  volume = {255},
  ISSN = {0377-2217},
  DOI = {10.1016/j.ejor.2016.05.048},
  number = {2},
  journal = {European Journal of Operational Research},
  publisher = {Elsevier BV},
  author = {Moshref-Javadi,  Mohammad and Lee,  Seokcheon},
  year = {2016},
  month = dec,
  pages = {604–619}
}

@article{Jozefowiez2014,
  title = {A branch-and-price algorithm for the multivehicle covering tour problem},
  volume = {64},
  ISSN = {0028-3045},
  DOI = {10.1002/net.21564},
  number = {3},
  journal = {Networks},
  publisher = {Wiley},
  author = {Jozefowiez,  Nicolas},
  year = {2014},
  month = sep,
  pages = {160–168}
}

@article{Glize2020,
  title = {Exact methods for mono-objective and Bi-Objective Multi-Vehicle Covering Tour Problems},
  volume = {283},
  ISSN = {0377-2217},
  DOI = {10.1016/j.ejor.2019.11.045},
  number = {3},
  journal = {European Journal of Operational Research},
  publisher = {Elsevier BV},
  author = {Glize,  Estèle and Roberti,  Roberto and Jozefowiez,  Nicolas and Ngueveu,  Sandra Ulrich},
  year = {2020},
  month = jun,
  pages = {812–824}
}

@article{Ropke2006,
  title = {An Adaptive Large Neighborhood Search Heuristic for the Pickup and Delivery Problem with Time Windows},
  volume = {40},
  ISSN = {1526-5447},
  DOI = {10.1287/trsc.1050.0135},
  number = {4},
  journal = {Transportation Science},
  publisher = {Institute for Operations Research and the Management Sciences (INFORMS)},
  author = {Ropke,  Stefan and Pisinger,  David},
  year = {2006},
  month = nov,
  pages = {455–472}
}

@article{Brimberg2023,
  title = {Variable Neighborhood Search: The power of change and simplicity},
  volume = {155},
  ISSN = {0305-0548},
  DOI = {10.1016/j.cor.2023.106221},
  journal = {Computers \& Operations Research},
  publisher = {Elsevier BV},
  author = {Brimberg,  Jack and Salhi,  Said and Todosijević,  Raca and Urošević,  Dragan},
  year = {2023},
  month = jul,
  pages = {106221}
}

@article{Hachicha2000,
  title = {Heuristics for the multi-vehicle covering tour problem},
  volume = {27},
  ISSN = {0305-0548},
  DOI = {10.1016/s0305-0548(99)00006-4},
  number = {1},
  journal = {Computers \& Operations Research},
  publisher = {Elsevier BV},
  author = {Hachicha,  Mondher and John Hodgson,  M and Laporte,  Gilbert and Semet,  Frédéric},
  year = {2000},
  month = jan,
  pages = {29–42}
}

@article{Gendreau1997,
  title = {The Covering Tour Problem},
  volume = {45},
  ISSN = {1526-5463},
  DOI = {10.1287/opre.45.4.568},
  number = {4},
  journal = {Operations Research},
  publisher = {Institute for Operations Research and the Management Sciences (INFORMS)},
  author = {Gendreau,  Michel and Laporte,  Gilbert and Semet,  Frédéric},
  year = {1997},
  month = aug,
  pages = {568–576}
}

@article{MORADI2024110730,
title = {Set Covering Routing Problems: A review and classification scheme},
journal = {Computers \& Industrial Engineering},
publisher = {Elsevier},
volume = {198},
pages = {110730},
year = {2024},
DOI = {10.1016/j.cie.2024.110730},
author = {Nima Moradi and Fereshteh Mafakheri and Chun Wang},
}

@article{allahyari2015hybrid,
  title = {A hybrid metaheuristic algorithm for the multi-depot covering tour vehicle routing problem},
  volume = {242},
  ISSN = {0377-2217},
  DOI = {10.1016/j.ejor.2014.10.048},
  number = {3},
  journal = {European Journal of Operational Research},
  publisher = {Elsevier BV},
  author = {Allahyari,  Somayeh and Salari,  Majid and Vigo,  Daniele},
  year = {2015},
  month = may,
  pages = {756–768}
}

@article{nedjati2017bi,
  title={Bi-objective covering tour location routing problem with replenishment at intermediate depots: Formulation and meta-heuristics},
  author={Nedjati, Arman and Izbirak, Gokhan and Arkat, Jamal},
  journal={Computers \& Industrial Engineering},
  volume={110},
  pages={191--206},
  year={2017},
  DOI={10.1016/j.cie.2017.06.004},
  publisher={Elsevier},
}

@article{groer2010library,
  title={A library of local search heuristics for the vehicle routing problem},
  author={Gro{\"e}r, Chris and Golden, Bruce and Wasil, Edward},
  journal={Mathematical Programming Computation},
  volume={2},
  number={2},
  pages={79--101},
  year={2010},
  DOI={10.1007/s12532-010-0013-5},
  publisher={Springer}
}

@article{funke2005local,
  title={Local search for vehicle routing and scheduling problems: Review and conceptual integration},
  author={Funke, Birger and Gr{\"u}nert, Tore and Irnich, Stefan},
  journal={Journal of heuristics},
  volume={11},
  number={4},
  pages={267--306},
  year={2005},
  DOI={10.1007/s10732-005-1997-2},
  publisher={Springer}
}

@article{NajiAzimi2012,
  title = {An Integer Linear Programming based heuristic for the Capacitated m-Ring-Star Problem},
  volume = {217},
  ISSN = {0377-2217},
  DOI = {10.1016/j.ejor.2011.08.026},
  number = {1},
  journal = {European Journal of Operational Research},
  publisher = {Elsevier BV},
  author = {Naji-Azimi,  Zahra and Salari,  Majid and Toth,  Paolo},
  year = {2012},
  month = feb,
  pages = {17–25}
}

@article{janinhoff2024out,
title = {Out-of-home delivery in last-mile logistics: A review},
journal = {Computers \& Operations Research},
volume = {168},
pages = {106686},
year = {2024},
DOI={10.1016/j.cor.2024.1066860},
author = {Lukas Janinhoff and Robert Klein and Daniela Sailer and Jim Morten Schoppa},
}

@article{lindauer-jmlr22a,
       author  = {Marius Lindauer and Katharina Eggensperger and Matthias Feurer and André Biedenkapp and Difan Deng and Carolin Benjamins and Tim Ruhkopf and René Sass and Frank Hutter},
       title   = {SMAC3: A Versatile Bayesian Optimization Package for Hyperparameter Optimization},
       journal = {Journal of Machine Learning Research},
       year    = {2022},
       volume  = {23},
       number  = {54},
       pages   = {1--9},
       }

@article{RANCOURT201568,
title = {Tactical network planning for food aid distribution in Kenya},
journal = {Computers \& Operations Research},
volume = {56},
pages = {68-83},
year = {2015},
issn = {0305-0548},
doi = {https://doi.org/10.1016/j.cor.2014.10.018},
author = {Rancourt, Marie-Eve and Cordeau, Jean-Francois and Laporte, Gilbert and Watkins, Ben},
}

@misc{WFP2024UNHAS,
  author       = {UNHAS},
  title        = {UNHAS Annual Review 2024},
  howpublished = {\url{https://www.wfp.org/publications/unhas-annual-review}},
  year         = {2025},
  note         = {Accessed: 23 April 2026},
}

@misc{heli_datasheet,
  author       = {BlueSkyRotor},
  title        = {Mil Mi-8AMT Data Sheet},
  howpublished = {\url{https://www.blueskyrotor.com/performance/datasheet/Mil/Mi_8-Mi_8-AMT}},
  year         = {2026},
  note         = {Accessed: 23 April 2026},
}

@misc{affected_individuals,
  author       = {OCHA},
  title        = {Afghanistan: Natural Disaster Incidents},
  howpublished = {\url{https://data.humdata.org/dataset/afghanistan-natural-disaster-incidents}},
  year         = {2025},
  note         = {Accessed: 23 April 2026},
}

@misc{population_data,
  author       = {WorldPop},
  title        = {Afghanistan - Spatial Distribution of Population (2015-2030)},
  howpublished = {\url{https://data.humdata.org/dataset/worldpop-population-counts-2015-2030-afg}},
  year         = {2025},
  note         = {Accessed: 23 April 2026},
}

@misc{min_calorie_intake,
  author       = {{Afghanistan Food Security Cluster}},
  title        = {Guideline on Food Security and Agriculture Cluster Response Packages},
  howpublished = {\url{https://fscluster.org/ethiopiafsc/document/guideline-food-cluster-response-packages}},
  year         = {2026},
  note         = {Accessed: 23 April 2026},
}

@misc{wfp_heb_response,
  author       = {WFP},
  title        = {Afghanistan WFP Flash Appeal: Floods - May 2024},
  howpublished = {\url{https://reliefweb.int/report/afghanistan/afghanistan-wfp-flash-appeal-floods-may-2024}},
  year         = {2025},
  note         = {Accessed: 23 April 2026},
}

@Article{tandem_elevation_data,
AUTHOR = {González, Carolina and Bachmann, Markus and Bueso-Bello, José-Luis and Rizzoli, Paola and Zink, Manfred},
TITLE = {A Fully Automatic Algorithm for Editing the TanDEM-X Global DEM},
JOURNAL = {Remote Sensing},
VOLUME = {12},
YEAR = {2020},
NUMBER = {23},
ARTICLE-NUMBER = {3961},
ISSN = {2072-4292},
DOI = {10.3390/rs12233961}
}

@Article{esa_terrain_data,
  author       = {Zanaga, Daniele and
                  Van De Kerchove, Ruben and
                  Daems, Dirk and
                  De Keersmaecker, Wanda and
                  Brockmann, Carsten and
                  Kirches, Grit and
                  Wevers, Jan and
                  Cartus, Oliver and
                  Santoro, Maurizio and
                  Fritz, Steffen and
                  Lesiv, Myroslava and
                  Herold, Martin and
                  Tsendbazar, Nandin-Erdene and
                  Xu, Panpan and
                  Ramoino, Fabrizio and
                  Arino, Olivier},
  title        = {ESA WorldCover 10 m 2021 v200},
  month        = "October",
  year         = "2022",
  publisher    = {Zenodo},
  version      = {v200},
  doi          = {10.5281/zenodo.7254221},
}

@misc{road_networks,
  author       = {HOT},
  title        = {Afghanistan Roads (OSM Export)},
  howpublished = {\url{https://data.humdata.org/dataset/hotosm_afg_roads}},
  year         = {2026},
  note         = {Accessed: 23 April 2026},
}

@misc{wfp_funding_shortfall,
  author       = {WFP},
  title        = {Food Security Impact of Reduction in WFP Funding},
  howpublished = {\url{https://www.wfp.org/publications/food-security-impact-reduction-wfp-funding}},
  year         = {2025},
  note         = {Accessed: 23 April 2026},
}

@misc{awd_cholera,
  author       = {ACAPS},
  title        = {Afghanistan: impact of flooding},
  howpublished = {\url{https://www.acaps.org/en/countries/archives/detail/afghanistan-impact-of-flooding}},
  year         = {2024},
  note         = {Accessed: 23 April 2026},
}

@misc{inaccessible_source_who,
  author       = {WHO},
  title        = {Afghanistan Flooding Situation Report No. 3 (14 May 2024)},
  howpublished = {\url{https://reliefweb.int/report/afghanistan/afghanistan-flooding-situation-report-no-3-14-may-2024}},
  year         = {2025},
  note         = {Accessed: 23 April 2026},
}

@misc{inaccessible_source_unicef,
  author       = {UNICEF},
  title        = {Flash Floods in Northern Afghanistan - UNICEF Humanitarian Situation Update},
  howpublished = {\url{https://www.unicef.org/afghanistan/documents/flash-floods-northern-afghanistan-unicef-humanitarian-situation-update-1}},
  year         = {2025},
  note         = {Accessed: 23 April 2026},
}

@misc{unicef_extreme_weather,
    author       = {OCHA}, 
    title = {Afghanistan: The alarming effects of climate change},
    howpublished = {\url{https://www.unocha.org/news/afghanistan-alarming-effects-climate-change}},
    year = {2023},
    note = {Accessed: 23 April 2026},
}

@misc{afg_household_data,
    author = {NSIA},
    title = {Estimated Population of Afghanistan 2025-26},
    howpublished = {\url{https://nsia.gov.af/library}},
    year = {2025},
    note = {Accessed: 23 April 2026},
}

@misc{UNOCHA2024GHO2025-july,
  author       = {OCHA},
  title        = {Global Humanitarian Overview 2025 - July Update},
  year         = {2025},
  institution  = {United Nations Office for the Coordination of Humanitarian Affairs},
  howpublished = {\url{https://reliefweb.int/report/world/global-humanitarian-overview-2025-july-update-snapshot-31-july-2025-enar}},
  note = {Accessed: 23 April 2026}
}

@misc{UNOCHA2024GHO2025,
  author       = {OCHA},
  title        = {Global Humanitarian Overview 2025},
  year         = {2025},
  institution  = {United Nations Office for the Coordination of Humanitarian Affairs},
  howpublished = {\url{https://www.unocha.org/publications/report/world/global-humanitarian-overview-2025-enarfres}},
  note = {Accessed: 23 April 2026}
}

@article{Kammoun2016,
  title = {An integration of mixed VND and VNS: the case of the multivehicle covering tour problem},
  volume = {24},
  ISSN = {1475-3995},
  DOI = {10.1111/itor.12355},
  number = {3},
  journal = {International Transactions in Operational Research},
  publisher = {Wiley},
  author = {Kammoun,  Manel and Derbel,  Houda and Ratli,  Mostapha and Jarboui,  Bassem},
  year = {2016},
  month = oct,
  pages = {663–679}
}

@misc{heb_density,
  author       = {WFP},
  title        = {What WFP Delivers: High-Energy Biscuits},
  howpublished = {\url{https://www.wfpusa.org/news/what-wfp-delivers-high-energy-biscuits/}},
  year         = {2025},
  note         = {Accessed: 23 April 2026},
}

@misc{cpu_benchmark,
  author       = {{CPU Benchmarks}},
  title        = {CPU Benchmarks: Intel i7-8700K vs i7-4790},
  howpublished = {\url{https://www.cpubenchmark.net/compare/3098vs2226/Intel-i7-8700K-vs-Intel-i7-4790}},
  year         = {2026},
  note         = {Accessed: 23 April 2026},
}
